\documentclass[11pt,a4paper]{article}
\usepackage{amssymb,fullpage,amsmath,mathrsfs}
\usepackage{pgfplots}
\usepackage{amsfonts,amsbsy,graphicx}
\usepackage{authblk,lineno}

\newtheorem{theorem}{Theorem}[section]

\newcommand{\eeref}[1]{(\ref{eqn:#1})}
\newcommand{\eelab}[1]{\label{eqn:#1}}
\newcommand{\ffref}[1]{\ref{fig:#1}}
\newcommand{\fflab}[1]{\label{fig:#1}}

\def\beq{\begin{equation}}
\def\eeq{\end{equation}}

\def\XXint#1#2#3{{\setbox0=\hbox{$#1{#2#3}{\int}$}
\vcenter{\hbox{$#2#3$}}\kern-.5\wd0}}

\def\strutdepth{\dp\strutbox}
\def \nw#1{\strut\vadjust{\kern-\strutdepth\vtop to0pt{\vss\hbox to\hsize {\hskip\hsize\hskip5pt$\leftarrow$\hss\strut}}}{\em \textcolor{blue}{#1}}}

\begin{document}

\title{The Evolution of Travelling Waves in a KPP Reaction-Diffusion Model with cut-off Reaction Rate. I. Permanent Form Travelling Waves.}

\author{A. D. O. Tisbury, D. J. Needham and A. Tzella\thanks{Address for correspondence:  Prof. D. J. Needham and Dr A. Tzella, School of Mathematics, University of Birmingham; email: d.j.needham@bham.ac.uk and a.tzella@bham.ac.uk}} 
\affil{School of Mathematics, University of Birmingham, Birmingham, B15 2TT, UK} 
\maketitle

\begin{abstract}
We consider  Kolmogorov--Petrovskii--Piscounov (KPP) type models in the presence of a discontinuous cut-off in reaction rate at concentration $u=u_c$. In Part I we examine permanent form travelling wave   solutions (a companion paper, Part II, is devoted to their evolution in the large time limit). For each fixed cut-off value $0<u_c<1$, we prove the existence of a unique permanent form travelling wave with a continuous and monotone decreasing propagation speed $v^*(u_c)$.  We extend previous asymptotic results  in the limit of small $u_c$ and present new asymptotic results in the limit of  large $u_c$ which are respectively obtained via the systematic use of matched and regular asymptotic expansions. The asymptotic results are confirmed against numerical results obtained for  the particular case of a cut-off Fisher reaction function. 
\end{abstract}

\section{Introduction}
\begin{figure}[t]
   \centering
   \begin{minipage}[b]{0.48\textwidth}
	     \includegraphics[width=\textwidth]{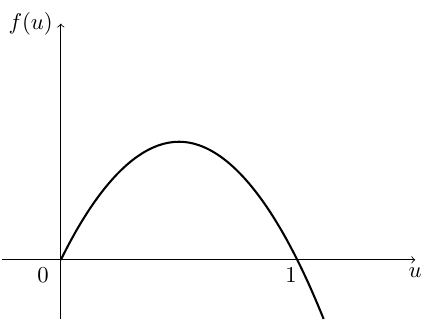}\\
	     \centering  (a)
   \end{minipage}
   \hfill
   \begin{minipage}[b]{0.49\textwidth}
	   \includegraphics[width=\textwidth]{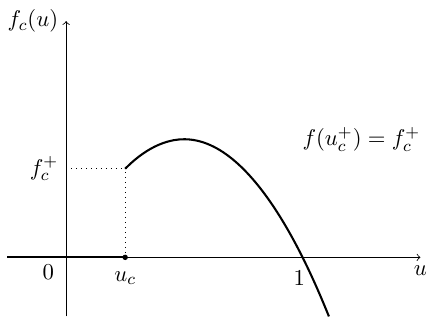}\\
  \centering  (b)
   \end{minipage}
 \caption{(a) A sketch of a  KPP-type reaction function.
 (b) A sketch of a cut-off KPP-type reaction function.}
 \label{fig:reaction_functions}
 \end{figure}

Travelling waves
arise in a wide range of applications in  
 mathematical chemistry and biology (for example, in combustion \cite{Williams1985}  and  in ecology, epidemiology and genetics  \cite{Fife1979,Murray}).
They describe the invasion of    chemical or biological reactions
and are usually established as a result of the interaction between molecular diffusion, local growth and saturation.
The simplest model that encapsulates this interaction is the
 Kolmogorov--Petrovskii--Piscounov (KPP) reaction--diffusion equation (also called Fisher-KPP    equation) \cite{Fisher1937,Kolmogorov_etal1937}.
 In one spatial dimension this describes the evolution of the concentration
 $u(x,t)$ as 
 \begin{linenomath}
 \begin{subequations}\eelab{KPP}
 \begin{align}
& u_t  = u_{xx} + f(u), \qquad (x,t) \in \mathbb{R} \times \mathbb{R}^+,\eelab{KPPa}\\
& u(x,0)=u_0(x), \eelab{KPPb} 
  \end{align}
where $u_0:\mathbb{R}\to\mathbb{R}$ is piecewise continuous
and smooth with limits $0$ and $1$ as $x\to\infty$ and $x\to-\infty$, respectively. This is typically supplemented
   with boundary conditions
   \beq\eelab{IBCb} %
     u(x,t)  \to  
	 \begin{cases}
		 1, & \text{as $x \to - \infty$}\\
		             0, & \text{as $x \to \infty$}
	\end{cases}				 
 \eeq
 \end{subequations}
 \end{linenomath}
  with these limits being uniform for  $t\in[0,T]$ and any $T>0$.
The  function
  $f: \mathbb{R} \to \mathbb{R}$ is a normalised KPP-type reaction function which satisfies  conditions that $f \in C^1(\mathbb{R})$ and 
  \begin{subequations} \eelab{KPPreaction}
 	\beq
	  	f(0) = f(1)= 0,\quad f'(0)  = 1, \quad f'(1)  <0
	\eeq	
	 	and in addition
		\beq
	  	 0<f(u)\leq u   \quad  \forall    u \in (0,1),\quad
	 	\quad f(u) < 0 \quad  \forall    u \in (1,\infty).
	  	\eeq
	  \end{subequations}  
	  A prototypical example of such a KPP reaction function is the Fisher reaction function \cite{Fisher1937} given by 
	  \begin{subequations}
	  \beq\eelab{Fisher}
	  	  f(u)=u(1-u).
	  \eeq
	  An illustration of $f(u)$ against $u$ is given in Figure \ref{fig:reaction_functions}(a). Another popular example of a KPP reaction function is
	   \beq\eelab{Fisher2}
	  f(u)=u(1-u^2).
	  \eeq
	  \end{subequations}

It is well-known \cite{AronsonWeinberger1975,Fife1979,Kolmogorov_etal1937,Smoller89} that the initial-boundary value problem  \eeref{KPP}  for the KPP equation  
supports a one-parameter family of  non-negative permanent form travelling wave   solutions of the form
    \beq\eelab{PTW}
    u(x,t)=U(y)=U(x-vt)  \quad \forall   (x,t) \in\mathbb{R}\times\mathbb{R}^+.
    \eeq
 These remain steady 
in time in a reference frame moving in the positive $x$ direction with speed
	$v\geq 0$ to be determined.
  Their existence and  uniqueness (up to linear translation in origin of the independent coordinate $y$) is established for
  \beq
  v\geq v_m=2,
  \eeq
  where $v_m$ denotes the minimum speed of propagation.
  This  is achieved by analysing the following nonlinear
 boundary value problem, namely,
 \begin{linenomath}
   \begin{subequations} \eelab{NBVP_KPPintro}
   	\begin{align}
   		  & U'' + vU' + f(U) = 0,\qquad-\infty< y <\infty,\\
		  & U(y)\geq 0,\qquad-\infty< y <\infty,\\
   		  & U(y)\to 
	   	 \begin{cases}
	   		 1, & \text{as $ y \to - \infty$}\\
	   		             0, & \text{as $y \to \infty$},
	   	\end{cases}
   	 \end{align}
   \end{subequations}
   \end{linenomath}
   where the dash denotes differentiation with respect to $y$. 
  This is  obtained by inserting  \eeref{PTW}
    into equation \eeref{KPPa} and using \eeref{IBCb}.
	The analysis is based on
	examining the existence of a unique heteroclinic orbit connecting the stable fixed
	point $(U,U')=(0,0)$ to the unstable fixed point $(U,U')=(1,0)$
	in the $(U,U')$ phase plane 
	of the equivalent two-dimensional dynamical system obtained from \eeref{NBVP_KPPintro}.
	It is also used to establish that
	$U(y)$  is monotone decreasing in  $y \in \mathbb{R}$. When translational invariance is fixed by requiring that $U(0)=1/2$, then  explicit expressions for the behaviour of the permanent form travelling wave  near the two fixed points are given by
	\begin{subequations}\eelab{PTW_KPP}
	\beq\eelab{PTW_KPPa}
	U(y)\sim
   	 \begin{cases}
   		 (A_\infty y+B_\infty)e^{-y}, & \text{as $ y \to \infty, v=v_m=2$}\\
   		             C_\infty e^{\alpha(v) y}, & \text{as $y \to \infty, v>v_m=2$}
   	\end{cases}
	\eeq
	and for all $v\geq v_m=2$,
	\beq\eelab{PTW_KPPb}
	 U(y)\sim 1-A_{-\infty} e^{\gamma(v) y},
	 \qquad  \textrm{as} \quad y \to -\infty,
	 \eeq
where
\beq
 \alpha(v)=\frac{1}{2}(-v+\sqrt{v^2-4})<0,
\qquad
\gamma(v)=1/2(-v+\sqrt{v^2+4|f'(1)|})>0,
 	 \eeq
	with
	 $A_\infty(>0)$, $B_\infty$, $C_\infty(>0)$  and $A_{-\infty}(>0)$
	being  
	globally determined constants,
	dependent on the nonlinearity of the boundary
	value problem \eeref{NBVP_KPPintro}
	 (see, for example, \cite{Fife1979,HadelerRothe1975}).
  \end{subequations}

	A key result is that the initial condition in \eeref{KPPb}
	determines  the permanent form travelling wave solution that emerges at large times.
	When $u_0(x)$ is
 sufficiently close to a Heaviside function, specifically, $u_0(x)\leq O(e^{-x})$ (meaning $O(e^{-x})$ or $o(e^{-x})$)  as $x\to\infty$,  
	the solution to the KPP initial-boundary value problem \eeref{KPP} converges at large times  to the
	permanent form travelling solution
	with minimum speed $v_m=2$
	 \cite{AronsonWeinberger1975,Kolmogorov_etal1937,Larson1978, Needham1992} at an algebraic
	 rate determined in \cite{Bramson1983,McKean1975,MerkinNeedham1993}.
 The mechanism which
 selects the speed of propagation  of the emerging permanent form travelling wave solution (as well as the
 rate of convergence)
 is based on the linearisation of the KPP equation \eeref{KPPa}
 at the leading edge of the travelling wave.
 There, the concentration $u$  is small and the dynamics are unstable.
As a result,  any modification of the dynamics
 near the leading edge of the travelling wave
 would
 invalidate this speed selection mechanism.
 
This is precisely  the case for the cut-off KPP model
that Brunet and Derrida \cite{BrunetDerrida1997} proposed and  considered.
Motivated by the discrete nature of chemical and biological phenomena at the microscopic level, 
they
took
 a reaction function that  
 is effectively deactivated at points where the concentration $u$ lies at or below a threshold value $u_c \in (0,1)$.  
This case corresponds to the cut-off  KPP equation   given by 
\begin{linenomath} 
  \begin{subequations} \eelab{cut-offKPP}
 	\begin{align}
   & u_t  = u_{xx} + f_c(u), \qquad (x,t) \in \mathbb{R} \times \mathbb{R}^+,\eelab{cut-offKPPa}\\
    & u(x,0)=u_0(x),
	\end{align}
     which is once more supplemented
    with the boundary conditions   
	\beq\eelab{cut-offIBCb}
     u(x,t)  \to 
	 \begin{cases}
	 			1, & \text{as $x \to - \infty$}\\
	             0, & \text{as $x \to \infty$}
	 		 \end{cases}
    \eeq
	uniformly for $t \in [0,T]$ for all $T >0$.
The main difference   is that the
 reaction function $f:\mathbb{R} \to \mathbb{R}$ in the KPP equation \eeref{KPP}
is replaced
with
  a cut-off reaction function
$f_c: \mathbb{R} \to \mathbb{R}$  given by
 \beq \eelab{BDreaction}
 f_c(u)=
 \begin{cases}
   			f(u), & \text{$u \in (u_c, \infty)$}\\
             0, & \text{$u \in (- \infty, u_c]$}
 		 \end{cases}
\eeq
where $f(u)$  satisfies the KPP conditions \eeref{KPPreaction}.
An illustration of $f_c(u)$ against $u$ is given in Figure \ref{fig:reaction_functions}b, with $f_c^+ = f_c(u_c^+)$ where $f_c(u_c^+)$ is the short notation  for $\lim_{u\to u_c^+} f_c(u)$. 
We remark that  $f_c(u)$ 
exhibits similarities with %
reaction functions arising in models of combustion 
in which  
$u_c$ represents 
an ignition temperature threshold \cite{Williams1985, Gordon2007}.
Focussing on the initial conditions
\beq
u_0(x)= 
\begin{cases}
			1, & \text{for $x < 0$}\\
            0, & \text{for $x \geq 0$}
		 \end{cases} 
\eeq
    \end{subequations}
    \end{linenomath}
we henceforth refer to this initial-boundary value problem as  IVP.
Brunet and Derrida \cite{BrunetDerrida1997} proposed   \eeref{cut-offKPP} 
as a   model of front 
  propagation arising in discrete systems of interacting 
particles. 
Such systems are for example, 
lattice models with discrete particles which make diffusive hops to neighbouring sites, and which have some birth-death  type of reaction \cite{Liggett2005}. 
In the continuum limit, 
obtained by allowing an arbitrarily large number of particles  per lattice site, 
Brunet and Derrida \cite{BrunetDerrida1997} conjectured that discreteness in concentration values can be represented by 
an effective  cut-off   $u_c$  where  $u_c$ may be viewed as  the effective mass of a single particle. 
The idea is that for $u<u_c$, 
diffusion dominates over growth.
Although the connection between \eeref{cut-offKPP}  and discrete systems of interacting 
particles is phenomenological,  model  \eeref{cut-offKPP} remains   useful in providing insight into their behaviour. 
Analysing the specific example \eeref{Fisher2}, 
Brunet and Derrida \cite{BrunetDerrida1997}
considered the behaviour of permanent form travelling wave solutions for small values of $u_c$. 
Their main result is a prediction for the propagation speed $v^*(u_c)$ of the \emph{unique}
permanent form travelling wave given by
\beq \eelab{speed_cutoff_BD}
v^*(u_c) \approx 2 - \frac{\pi^2}{(\ln u_c)^2},  \quad \textrm{as} \quad  u_c\to 0^+,
 \eeq
which they obtained using
a  two-region informal point patching procedure  (see also \cite{Kessler_etal1998} where \eeref{speed_cutoff_BD} was compared against numerical simulations of lattice particle models).
This   significant result demonstrates the strong influence of
a cut-off
on the value of $v^*(u_c)$ for small values of $u_c$.
The same approximation to $v^*(u_c)$  
has also been obtained via an alternative variational approach in
 \cite{BenguriaDepassier2007,Benguria_etal2008}. 
{Subsequently,    a more rigorous approach was employed by Dumortier, Popovic and Kaper  \cite{Dumortier_etal2007}  
who used geometric desingularisation, 
   to prove the existence and uniqueness  of a permanent form travelling wave
with }
\beq\eelab{speed_cut-off_Dumortier_etal}
v^*(u_c) \sim 2- \frac{\pi^2}{(\ln u_c)^2}+O\left(\frac{1}{|\ln u_c|^3}
\right),
\quad \textrm{as} \quad  u_c\to 0^+.
\eeq
{All these results have restricted validity to the small $u_c$ limit  
 with specific choices of cut-off KPP-type reaction function  \eeref{BDreaction},   the most common  based on $f(u)$ given by 
\eeref{Fisher2}}\footnote{There have been a number of results obtained for  other  cut--off reaction functions    (see for example, \cite{Gordon2007,Dumortier_etal2010,Popovic2011,DumortierKaper2012}),  but we focus on the  cut-off KPP-type  reaction functions.}. 
Expression \eeref{speed_cut-off_Dumortier_etal} was found in \cite{Dumortier_etal2007}  to be generic
when considering a slightly more general 
 class of cut-off  KPP-type reaction functions, namely 
  identical to (8d) when $u\in(u_c,\infty)$ but has $f_c(u)=o(1)$
uniformly for $u\in[0,u_c]$ as $u_c\to 0^+$.

There are a number of fundamental questions that remain.
The first question concerns the existence and uniqueness of a permanent form travelling wave solution for \emph{arbitrary} threshold values $u_c$ and KPP reaction functions $f(u)$.
The second question concerns the propagation speed of such permanent form travelling wave  solutions for arbitrary threshold values $u_c$. The third question is with regard to the shape of the permanent form travelling wave solution. 
 The fourth question concerns 
a systematic approach that captures the leading as well as higher order corrections to the asymptotic 
behaviour of the speed  and shape of the permanent form travelling wave solution as $u_c\to 0^+$   and $u_c\to 1^-$.   
The second limit may be less relevant  for discrete systems of interacting particles. 
It is however relevant in models of combustion  since the ignition 
temperature that determines the cut-off is not necessarily small \cite{Williams1985}.
A final question concerns the evolution in time to the permanent form travelling wave solution via the initial boundary value problem IVP.
 Part I of this series of papers addresses the first four of these questions while part II addresses the fifth and last question. 
In particular, we study
classical solutions $u: \mathbb{R} \times \overline{\mathbb{R}}^+ \to \mathbb{R}$ to IVP
for the   cut-off KPP  equation \eeref{cut-offKPP}.
In this paper we proceed as follows.
In Section 2,
we re-formulate IVP  as a moving boundary problem. We then make a simple coordinate transformation to consider an equivalent initial-boundary value problem that we refer to as QIVP.
In section 3, we examine the possibility that QIVP supports
permanent form travelling wave solutions
$U_T(y)=U_T(x-vt)$ where $U_T \in C^1(\mathbb{R}) \cap C^2(\mathbb{R} \setminus \{ 0 \})$
 satisfies the nonlinear boundary value problem,
 \begin{linenomath}
\begin{subequations} \eelab{BVP_cut-offintro}
\begin{align} 
& U_T''  + vU_T' + f_c(U_T) = 0,  \qquad y \in \mathbb{R} \setminus \{0 \},  \\
& U_T \geq u_c \quad \forall y < 0, \qquad  0 \leq U_T \leq u_c, \quad \forall  y > 0, \eelab{NBVP_KPPaintro} \\
&  U_T(0) = u_c, \eelab{NBVP_KPPbintro}  \\
& U_T(y) \to
 \begin{cases}
 			1, & \text{as $y \to - \infty$}\\
             0, & \text{as $y \to \infty$}
 		 \end{cases}
 \end{align}
\end{subequations}
\end{linenomath}
We establish the following theorem.
\begin{theorem} \label{theorem}
For each fixed $u_c \in (0,1)$, QIVP has a unique permanent form travelling wave solution $U_T:\mathbb{R} \to \mathbb{R}$, with the propagation speed given by $v^*(u_c)$. Here $v^*:(0,1) \to \mathbb{R}^+ $ is continuous and monotone decreasing, with
\[
v^*(u_c) \to
\begin{cases}
			0, & \text{as $u_c \to 1^-$}\\
            2, & \text{as $u_c \to 0^+$}
		 \end{cases}
\]
 where $2$ is the minimum propagation speed
 of the permanent form travelling wave solution
 in the absence of cut-off ($u_c =0$).
In addition, $U_T(y)$ is strictly monotone decreasing for $y \in \mathbb{R}$, with $U_T(0)=u_c$, and
\begin{linenomath}
\begin{subequations}
\begin{align}
&  U_T''(0^+)-U_T''(0^-) = - f_c^+ ,   \\
&  U_T(y) = u_c e^{- v^*(u_c)y} \quad \forall  y \in \overline{\mathbb{R}}^+,    \eelab{ypos} \\
&  U_T(y) \sim 1 - A_{-\infty} e^{\lambda_+ ( v^*(u_c))y} \quad \textrm{as} \quad y \to - \infty    ,
\end{align}
\end{subequations}
\end{linenomath}
for some global constant $A_{-\infty}>0$ (which depends upon $u_c$), and
\[
\lambda_+(v) =  \frac{1}{2}
\left( - v + \sqrt{v^2 + 4 |f_c'(1)| }
\right)
>0.
\]
Furthermore,
\beq\eelab{speed_cut-off_largeuc}
v^*(u_c) \sim |f_c'(1)|^{\frac{1}{2}} (1 - u_c) \qquad \textrm{as} \quad u_c \to 1^- .
\eeq
 \end{theorem}

In sections 4 and 5 we use matched asymptotic expansions to develop the detailed asymptotic structure to the permanent form travelling wave
solutions as $u_c \to 0^+$ and as $u_c \to 1^-$ respectively.
These are used to obtain
 higher-order corrections to \eeref{speed_cut-off_Dumortier_etal} and \eeref{speed_cut-off_largeuc} in a systematic manner.
In the first limit, the analysis
is carried out on the direct problem (rather than the phase plane).
It highlights that
higher-order corrections
are controlled by  two global constants $A_\infty$ and $B_\infty$
associated with the minimum speed of permanent form travelling wave
solution to the non cut-off KPP problem \eeref{KPP}.
These global constants  represent the nonlinearity in the problem
when $u_c$ is small.
The analysis is readily generalised to degenerate
and singular  KPP conditions, obtained for example
when $f'(0)=0$ or
 $f(u)\sim u^{1/2}$ as $u\to 0^+$, respectively.
Section 6 presents numerical examples for the specific Fisher  cut-off reaction function.
The paper concludes with a discussion in Section 7.
\section{Formulation of Evolution Problem QIVP}
\label{sec:QIVP}
\begin{figure}
\centerline{\includegraphics[width=0.61\linewidth]{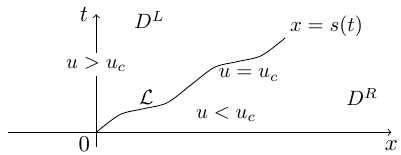}}
\caption{A sketch of the evolution of the boundary $u=u_c$ in the moving boundary problem.}
\label{fig:MIVP}
\end{figure}

Due to the discontinuity in $f_c(u)$ at $u=u_c$, it is convenient to re-structure  IVP as a moving boundary problem. To this end, we introduce the domains:
\begin{linenomath}
\begin{subequations} \eelab{domain}
\begin{align}
	 D^L = \{ (x,t) \in \mathbb{R} \times \mathbb{R}^+:  x < s(t)\} ,   \\
 D^R = \{ (x,t) \in \mathbb{R} \times \mathbb{R}^+:  x > s(t)\} ,
\end{align}
and the curve
\beq
\quad \mathcal{L} = \{ (x,t) \in \mathbb{R} \times \mathbb{R}^+: x = s(t)  \}    ,
\eeq
\end{subequations}
\end{linenomath}
that describes the moving boundary between the two domains. The boundary is expressed in terms of $s(t)$ which satisfies  $u(s(t),t)=u_c$, with $u \geq u_c$ in
$\overline{D}^L$
and $u \leq u_c$ in $\overline{D}^R$.
In this context, a classical solution will have $u: \mathbb{R} \times \overline{\mathbb{R}}^+ \to \mathbb{R}$ and $s:\overline{\mathbb{R}}^+ \to \mathbb{R}$
such that,
\begin{linenomath}
\begin{subequations}
\begin{align}
& u \in C ( \mathbb{R} \times \overline{\mathbb{R}}^+  \setminus \{ (0,0) \} ) \cap  C^{1,1} ( \mathbb{R} \times \mathbb{R}^+ ) \cap C^{2,1}  ( D^L \cup D^R),  \eelab{eq: continuity condition on u in moving boundary problem chap 2}  \\
& s \in  C^1  ( \mathbb{R}^+ ),    \eelab{eq: continuity condition on s in moving boundary problem chap 2}   \\
& s(0^+)  = 0.    \eelab{sinitial}
\end{align}
\end{subequations}
\end{linenomath}
The moving boundary problem is then formulated as follows,
\begin{linenomath}
\begin{subequations} \eelab{MBP}
 \begin{align}
& u_t  = u_{xx} + f_c(u), \qquad (x,t) \in D^L \cup D^R , \\
&  u  \geq u_c\;\textrm{in   $\overline{D}^L$}, \qquad u \leq u_c \;\textrm{in $\overline{D}^R$},      \\
& u(x,0)  =
\begin{cases}
			1, & \text{for $x < 0$}\\
            0, & \text{for $x\geq 0$}
		 \end{cases}
\\
& u(x,t)=
\begin{cases}
			1, & \text{as $x \to - \infty$}\\
            0, & \text{as $x\to\infty$}
		 \end{cases}
 \end{align}
 uniformly for $t \in [0,T]$ for all  $T >0$ and
 \begin{align}
& u(s(t),t)  = u_c, \qquad t \in \mathbb{R}^+, \\
& u_x(s(t)^+,t) = u_x(s(t)^-,t), \qquad t \in \mathbb{R}^+, \\
&  s(0^+)  = 0.
\end{align}
\end{subequations}
\end{linenomath}
The situation is illustrated in Figure \ref{fig:MIVP}.
It is now convenient to make the simple coordinate transformation $(x,t) \to (y,t)$ with $y=x-s(t)$. We then introduce the following domains:
\beq
Q^L = \mathbb{R}^- \times \mathbb{R}^+, \qquad Q^R = \mathbb{R}^+ \times \mathbb{R}^+,
\eeq
with $u: \mathbb{R} \times \overline{\mathbb{R}}^+ \to \mathbb{R}$ and $s:\overline{\mathbb{R}}^+ \to \mathbb{R}$ such that
\begin{linenomath}
\begin{subequations} \eelab{u_full_conditions}
\begin{align}
& u \in C  (  \mathbb{R} \times \overline{\mathbb{R}}^+  \setminus \{ (0,0) \} ) \cap C^{1,1}  ( \mathbb{R} \times \mathbb{R}^+) \cap C^{2,1}  ( Q^L \cup Q^R ),  \\
& s \in  C^1  ( \mathbb{R}^+ ).    \eelab{eq: continuity condition on s after transformation chap 2}
\end{align}
\end{subequations}
\end{linenomath}
The equivalent problem
to \eeref{MBP} is then given by
\begin{linenomath}
\begin{subequations} \eelab{QIVP1}
\begin{align}
& u_t - \dot{s}(t) u_y  = u_{yy} + f_c(u), \quad (y,t) \in Q^L \cup Q^R ,    \eelab{QIVP1a}  \\
& u  \geq u_c  \;\textrm{in  $\overline{Q}^L$}, \quad u \leq u_c  \;\textrm{in  $\overline{Q}^R$},      \\
& u(y,0)  =
 \begin{cases}
 			1, & \text{$y< 0$}\\
             0, & \text{$y\geq 0$}
 		 \end{cases}, \eelab{QIVP1b} \\
& u(y,t)  \to
 \begin{cases}
 			1, & \text{as $y\to-\infty$}\\
             0, & \text{as $y\to\infty$} \end{cases}, \eelab{QIVP1f}
 \end{align}
uniformly for $t \in [0,T]$ for all $T >0$ and
\begin{align}
& u(0,t)  = u_c, \quad t \in \mathbb{R}^+,  \eelab{QIVP1c}  \\
& u_y(0^+,t)  = u_y(0^-,t), \quad t \in \mathbb{R}^+,  \eelab{QIVP1d} \\
& s(0^+)  = 0, \eelab{QIVP1e}
\end{align}
\end{subequations}
\end{linenomath}
where the dot denotes differentiation with respect to time, $t$.
This initial-boundary value problem will henceforth be referred to as QIVP. On using the classical maximum principle and comparison theorem (see, for example, \cite{Fife1979} and \cite{AronSer67}), together with translational invariance in $y$, and the regularity in \eeref{u_full_conditions}, we can establish the following basic qualitative properties concerning QIVP, namely,
\begin{linenomath}
\begin{subequations}
\begin{align}
& 0 < u(y,t) < u_c \quad \forall  (y,t) \in Q^R,   \label{R1 chap 2}  \\
& u_c < u(y,t) < 1 \quad \forall (y,t) \in Q^L,    \\
& u(y,t) \mbox{ is strictly monotone decreasing in } y \in \mathbb{R} \quad \forall  t \in \mathbb{R}^+.  \eelab{u_conditions}
\end{align}
In addition, via the partial differential equation \eeref{QIVP1a} and the regularity conditions \eeref{u_full_conditions}, we have
\begin{align}
& \lim_{y \to 0^+} u_{yy} (y,t)  = \lim_{y \to 0^+} \left( u_t(y,t) - \dot{s}(t) u_y(y,t) \right)  \eelab{u_conditions2} \\
& \qquad \qquad \qquad = - \dot{s}(t) u_y(0,t) \quad \forall  t \in \mathbb{R}^+ , \nonumber   \\
& \lim_{y \to 0^-} u_{yy} (y,t) = \lim_{y \to 0^-} \left( u_t(y,t) - \dot{s}(t) u_y(y,t) -f(u(y,t)) \right)  \eelab{u_conditions3} \\
& \qquad \qquad \qquad = - \dot{s}(t) u_y(0,t) - f_c^+ \quad \forall  t \in \mathbb{R}^+ , \nonumber
\end{align}
with the limits in \eeref{u_conditions2} and \eeref{u_conditions3} being uniform for $t \in [t_0, t_1]$ (for any $0 < t_0 <t_1$). It follows from \eeref{u_conditions2} and \eeref{u_conditions3} that
\begin{equation} \eelab{u_conditions4}
\left[ u_{yy} (y,t) \right]^{y = 0^+}_{y = 0^-} = f_c^+ \quad \forall  t \in \mathbb{R}^+ ,
\end{equation}
whilst, using \eeref{u_conditions}, \eeref{u_conditions3} and the regularity condition \eeref{u_full_conditions} we establish that
\begin{align}
u_y(y,t) < 0 \quad \forall  (y,t) \in \mathbb{R} \times \mathbb{R}^+ .
\end{align}
\end{subequations}
\end{linenomath}
The remainder of this
paper and its companion (part II)
concentrates on the analysis of QIVP. Specifically, in this paper we consider the existence and uniqueness of permanent form travelling wave solutions  to QIVP
including their asymptotic behaviour in the limits of $u_c \to 0^+$ and $u_c \to 1^-$ via the method of matched and regular asymptotic   expansions.

\section{Permanent Form Travelling Waves in QIVP}
We anticipate that
as $t \to \infty$, a permanent form travelling wave solution   will develop in the solution to QIVP, advancing with a (non-negative) propagation speed,
allowing for the transition between the fully reacted state, $u=1$ as $y \to - \infty$, to the unreacted state, $u=0$ as $y \to \infty$. Therefore, in this section we focus attention on the possibility of QIVP supporting permanent form travelling wave solutions (henceforth referred to as PTW solutions).
We begin by establishing the existence and uniqueness of a PTW to QIVP for each fixed $u_c \in (0,1)$,
denoting the unique propagation speed by $v=v^*(u_c)$.
We then  consider limiting values of $v^*(u_c)$ as $u_c\to 0^+$ and $u_c\to 1^-$.
The results  established in this section provide proof of
Theorem \ref{theorem} as stated in the introduction.

\subsection{The Existence and Uniqueness of a PTW Solution to QIVP}
\label{sec:existence}

A PTW solution to QIVP, with constant speed of propagation $v \geq 0$,
is a steady state solution to QIVP with $u: \mathbb{R} \times \mathbb{R}^+ \to \mathbb{R}$ and $s: \overline{\mathbb{R}}^+ \to \mathbb{R}$ such that
\begin{linenomath}
\begin{align}
& u(y,t) = U_T(y) \quad \forall   (y,t) \in \mathbb{R} \times \mathbb{R}^+,   \eelab{eq:travelling wave notation}  \\
& \dot{s}(t) = v \quad \forall  t \in \mathbb{R}^+,  \eelab{eq:relation between sdot and v}
\end{align}
\end{linenomath}
where $U_T \in C^1(\mathbb{R}) \cap C^2(\mathbb{R} \setminus \{ 0 \})$ and $v \geq 0$ satisfy the nonlinear boundary value problem,
%
\begin{linenomath}
\begin{subequations}\eelab{BVP_cutoff}
\begin{align}
& U_T''  + vU_T' + f_c(U_T) = 0,  \qquad y \in \mathbb{R} \setminus \{0 \},  \\
& U_T \geq u_c \quad \forall  y < 0, \qquad  0 \leq U_T \leq u_c \quad \forall  y > 0, \eelab{NBVP_KPPa} \\
& U_T(0) = u_c, \eelab{NBVP_KPPb}  \\
& U_T(y) \to
 \begin{cases}
 			1, & \text{for $y \to - \infty$}\\
             0, & \text{for $y \to \infty$}
 		 \end{cases} \eelab{farfield}
 \end{align}
\end{subequations}
\end{linenomath}
where the dash denotes differentiation with respect to $y$.
The nonlinear boundary value problem \eeref{BVP_cutoff} can be thought of as a nonlinear eigenvalue problem with the eigenvalue being the propagation speed $v \geq 0$.

\addtocounter{equation}{1}

It is convenient to consider the ordinary differential equation \eeref{BVP_cutoff} as the following equivalent autonomous
first-order
two-dimensional dynamical system, with $\alpha=U_T$ and $\beta=U_T'$, namely,
\[
  \begin{array}{c l}  \eelab{phaseplane} \tag{26}
    \alpha' & = \beta ,  \quad  \\
   \beta' & = - v \beta - f_c(\alpha) , \quad
  \end{array}  \]
We will analyse this dynamical system in the $(\alpha,\beta)$ phase plane for $v\geq 0$.
In particular,
it is straightforward to establish that the existence of a solution to \eeref{BVP_cutoff}
is equivalent to the existence of a heteroclinic connection in the $(\alpha,\beta)$ phase plane, for the dynamical system \eeref{phaseplane}, which connects the equilibrium point $(1,0)$, as $y \to - \infty$, to the equilibrium point $(0,0)$, as $y \to \infty$ (the translational invariance is then fixed by condition \eeref{NBVP_KPPb}
which requires that $\alpha(0)=u_c$).
From \eeref{NBVP_KPPa},
this heteroclinic connection must remain in the $\alpha \geq 0$ half plane of the $(\alpha,\beta)$ phase plane, which we denote by $R^+=\{ (\alpha,\beta) : (\alpha, \beta) \in \overline{\mathbb{R}}^+ \times \mathbb{R} \}$. We henceforth focus
on this region of the $(\alpha,\beta)$ phase plane.

\begin{figure}
\centerline{\includegraphics[width=0.55\linewidth]{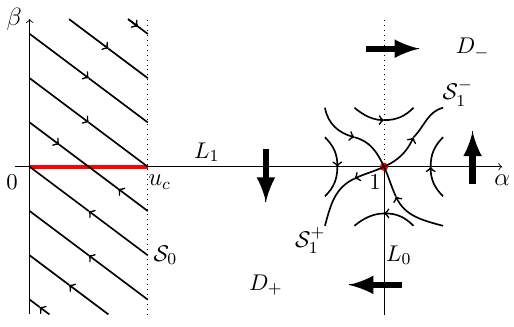}}
\caption{The local phase portrait for the equilibrium points of the dynamical system \eeref{phaseplane}. The thick black arrows denote the direction of the vector field $\vec{Q} (\alpha,\beta)$ along the line segments $L_0$, $L_1$ and on the boundary of $D_-$.
}
\label{fig:phase_portrait}
\end{figure}

However, before we proceed further, it is first worth considering the effect of introducing the cut-off into the reaction function on the dynamical system \eeref{phaseplane}. To that end, we introduce the function $\vec{Q} : \mathbb{R}^2 \to \mathbb{R}^2$ where $\vec{Q} (\alpha, \beta) $ is given by
\begin{equation} \eelab{vectorfield}
\vec{Q} (\alpha, \beta) = ( \beta , - v \beta - f_c (\alpha)) ,
\end{equation}
to represent the vector field generating the dynamical system \eeref{phaseplane}. We observe that, in the $(\alpha,\beta)$ phase plane, the effect of the discontinuity in $f_c(\alpha)$ across the line $\alpha = u_c$ is simply to $\it{refract}$ the phase paths passing through this line. In particular, for each $\beta \in \mathbb{R}$, there is exactly one phase path passing through $(u_c, \beta)$, which has tangent vectors, $\vec{Q} (u_c^-,\beta) = (\beta , -v \beta)$ and $\vec{Q} (u_c^+,\beta) = (\beta , -v \beta - f_c^+)$. Thus, the refraction vector for the phase paths which cross the line $\alpha = u_c$ is
\beq
\vec{R} (u_c, \beta)  = \vec{Q} (u_c^+,\beta) - \vec{Q} (u_c^-,\beta)  = ( 0 , -f_c^+). \eelab{refraction vector}
\eeq
We observe that the refraction vector \eeref{refraction vector} is independent of $(\beta, v) \in \mathbb{R} \times \overline{\mathbb{R}}^+$ and depends continuously on $u_c \in (0,1)$. It follows that
\begin{equation} \eelab{refraction_vector_limit}
\vec{R} (u_c, \beta) \to \vec{0} \quad \textrm{as} \quad u_c \to  0 ,
\end{equation}
uniformly in $(\beta, v) \in \mathbb{R} \times \overline{\mathbb{R}}^+$.
 After determining the effect of the discontinuity on the phase paths of the dynamical system \eeref{phaseplane} in $R^+$, we next consider the equilibrium points of \eeref{phaseplane} in $R^+$. These are readily found to be at locations
\beq	
 \vec{e_{a}} = (a,0) \quad \textrm{for each } a \in [0,u_c], \qquad \vec{e_1}=(1,0).
\eeq
 We begin by examining the local phase portrait in the neighbourhood of the equilibrium point $\vec{e_1}$. We find that $\vec{e_1}$ is a hyperbolic equilibrium point. Moreover, $\vec{e_1}$ is a saddle point with eigenvalues
\begin{equation}  \eelab{evals_e1}
\lambda_{\pm}(v) = \frac{1}{2} \left( - v \pm \sqrt{v^2 + 4 |f_c'(1)|} \right).
\end{equation}
The associated local 
 one-dimensional unstable and stable manifolds of 
$\vec{e_1}$ are, respectively, given by
\begin{equation} \eelab{unstable_manifold}
\beta(\alpha) = - \lambda_{\pm}(v)(1 - \alpha). 
\end{equation}
We denote the phase path which forms the part of the 
 (one-dimensional) 
 unstable  manifold entering $D_+= \{ (\alpha,\beta): 0 < \alpha < 1, \beta < 0 \}$ as $\mathcal{S}_{1}^{+}$. Similarly, we denote $\mathcal{S}_{1}^{-}$ as the phase path which forms part of the (one-dimensional)  unstable manifold entering $D_-= \{ (\alpha,\beta): \alpha > 1, \beta > 0  \}$. The situation is illustrated in Figure \ref{fig:phase_portrait}.
We next determine the local phase portrait of the equilibrium points $\vec{e_a}$ for each $a \in [0,u_c]$. For $ a \in (0,u_c)$ and $v>0$, each of the equilibrium points $\vec{e_a}$ is non-hyperbolic with a single (one-dimensional) stable manifold in $R^+$ given by $\{ (\alpha,\beta): \beta = -v(\alpha - a);\; 0 \leq \alpha \leq u_c \}$. Also, the equilibrium point $\vec{e_0}$ is non-hyperbolic with a single (one-dimensional) stable manifold in $\mathbb{R}^+$ which we will denote by
\begin{equation}  \eelab{stable manifold of e0}
\mathcal{S}_{0} = \{ (\alpha,\beta): \beta = -v \alpha  ; \; 0 \leq \alpha \leq u_c \} .
 \end{equation}
Finally, the equilibrium point $\vec{e_{u_c}}$ is again non-hyperbolic, and, for $0 \leq \alpha \leq u_c$, has a single  (one-dimensional)  stable manifold in $R^+$ given by
$\{(\alpha,\beta): \beta = -v(\alpha - u_c) ;\; 0 \leq \alpha \leq u_c \}$. In fact, the collection of phase paths of the dynamical system \eeref{phaseplane} in the region $\{ (\alpha,\beta):  0 \leq \alpha \leq u_c,\; \beta \leq 0 \}$ is given by the family of curves $\beta = c - v \alpha$, for each $c \in \mathbb{R}$. This is illustrated in Figure \ref{fig:phase_portrait}.
Next, for the line segment $\{ (\alpha,\beta): \alpha = 1,\; \beta > 0 \}$, we observe the following,
\begin{equation} \eelab{D-}
\vec{Q}(\alpha , \beta) \cdot (1,0) = \beta > 0 .
\end{equation}
Similarly, for the line segment $\{ (\alpha,\beta): \alpha > 1,\; \beta =0 \}$, we observe that
\begin{equation} \eelab{D-2}
\vec{Q}(\alpha , \beta) \cdot (0,1) = -f_c(\alpha) > 0 .
\end{equation}
Together with the local structure at the equilibrium point $\vec{e_1}$, we conclude from \eeref{D-} and \eeref{D-2} that the region $D_-$ is a strictly positively invariant region for the dynamical system \eeref{phaseplane}.
We now examine the line segments $L_0=\{(\alpha,\beta): \alpha = 1, \; \beta < 0\}$ and $L_1=\{(\alpha,\beta):  u_c < \alpha < 1,\; \beta =0 \}$, we observe that
\beq
\vec{Q}(\alpha , \beta) \cdot (-1,0) = - \beta > 0 \quad  \forall  (\alpha,\beta) \in L_0 , 
\quad
 \vec{Q}(\alpha , \beta) \cdot (0,-1) =  f_c(\alpha) > 0   \quad  \forall  (\alpha,\beta) \in L_1. \eelab{D+2}
\eeq
In addition, for $v>0$, we observe that for all $(\alpha,\beta) \in R^+$
\begin{equation} \eelab{bendixson}
\nabla \cdot \vec{Q}(\alpha,\beta) = - v < 0.
\end{equation}
Thus, for any $v>0$, it follows from the Bendixson negative criterion (see, for example, \cite{Verhulst90}) that \eeref{phaseplane} has no periodic orbits, homoclinic orbits or heteroclinic cycles in $R^+$. Finally, we observe that at each $(\alpha, \beta) \in R^+ \setminus ( \{\vec{e_1} \} \cup \{\vec{e_a}: 0 \leq a \leq u_c \})$ the vector field $\vec{Q}(\alpha, \beta)$ rotates continuously clockwise for increasing $v \geq 0$. At the equilibrium point $\vec{e_1}$, the unstable manifold $\mathcal{S}_{1}^{+}$ rotates clockwise for increasing $v \geq 0$,
as does the stable manifold $\mathcal{S}_{0}$ at the equilibrium point $\vec{e_0}$.
As the phase path $\mathcal{S}_1^-$ enters $D_-$ on leaving $\vec{e_1}$, and we have established that $D_-$ is a strictly positively invariant region for the dynamical system \eeref{phaseplane}, we conclude that this cannot correspond to a heteroclinic connection between $\vec{e_1}$ and $\vec{e_0}$. Thus, at any $v \geq 0$, the existence of a heteroclinic connection in $R^+$ connecting $\vec{e_1}$, as $y \to - \infty$, to $\vec{e_0}$, as $ y \to \infty$, is equivalent to the phase path $\mathcal{S}_1^+$, leaving $\vec{e_1}$, being coincident with the phase path $\mathcal{S}_0$, entering $\vec{e_0}$. It also follows that, at those $v \geq 0$ when such a heteroclinic connection exists, then it is unique.

We are now in a position to investigate for which values of $v \geq 0$, if any, the required heteroclinic connection exists in $R^+$. When $v=0$, it follows directly from \eeref{phaseplane} that the phase path $\mathcal{S}_1^{+}$ has graph $(\alpha,\beta_{0}(\alpha))$ where
\begin{equation} \eelab{beta_path}
\beta_{0}(\alpha) = - \left( 2 \int^1_{\alpha}  f_c(\gamma) d \gamma   \right)^{\frac{1}{2}},
\end{equation}
for $\alpha \in [0,1]$. Thus, $\beta_{0}(\alpha)$ is (non-positive) non-decreasing for $\alpha \in [0,1]$ with
\begin{equation}\eelab{initial_condition_beta0}
\beta_{0}(0) =  - \left( 2 \int^1_{u_c}  f_c(\gamma) d \gamma   \right)^{\frac{1}{2}} < 0
\quad
\text{and}
\quad
\beta_{0}'(1)  =  \left( -f_c'(1)    \right)^{\frac{1}{2}}.
\end{equation}
We also note that $\beta_0(\alpha)$ is continuous and differentiable except for a jump in derivative at $\alpha=u_c$ when $\beta_0'(u_c^+)=-f_c^+/\beta_0(0)$ whilst $\beta_0'(u_c^-)=0$.

\begin{figure}
\centerline{\includegraphics[width=0.6\linewidth]{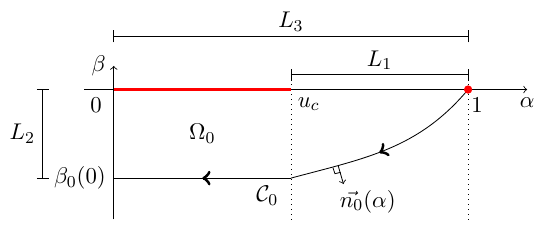}}
\caption{The phase path $\mathcal{C}_0$ which forms part of the unstable manifold of the equilibrium point $\vec{e_1}=(1,0)$ of the dynamical system \eeref{phaseplane} when $v=0$.}
\label{fig:beta0}
\end{figure}

We denote the phase path $\mathcal{S}_1^+|_{v=0}$ as $\mathcal{C}_0$, and note from \eeref{beta_path} that $\mathcal{C}_0 \subset \overline{D}^+$ as illustrated in Figure \ref{fig:beta0}. We conclude from \eeref{initial_condition_beta0} that when $v=0$ no heteroclinic connection exists from $\vec{e_1}$ to $\vec{e_0}$.
Moreover, it follows from the rotational properties of the vector field $\vec{Q}(\alpha, \beta)$ with increasing $v \geq 0$,
as discussed earlier, that, for each $v>0$, we have
\begin{equation} \eelab{L0}
\vec{Q}(\alpha , \beta_0(\alpha)) \cdot \vec{{n}_0}(\alpha) < 0,
\end{equation}
for all $\alpha \in [0,1)$, where $\vec{{n}_0}(\alpha)$ is the unit normal to $\mathcal{C}_0$ for $\alpha\in(u_c,1]$  as shown in Figure
\ref{fig:beta0}.
We define the line segments $L_2 =  \{(\alpha,\beta):  \alpha = 0,\; \beta_0(0) < \beta < 0 \}$ and $L_3 =\{(\alpha,\beta):   0 \leq \alpha \leq 1,\; \beta = 0 \}$ and denote the region $\Omega_0 \subset D_+$ as that region bounded by $\partial \Omega_0 = L_2 \cup L_3 \cup \mathcal{C}_0 $. We observe, via the rotational properties of $\mathcal{S}_1^+$ at $\vec{e_1}$ with increasing $v \geq 0$, that for any $v > 0$,  
$\mathcal{S}_1^+|_{v}$ enters $\Omega_0$ on leaving $\vec{e_1}$. Moreover, we recall, $\overline{\Omega}_0$ contains no periodic orbits, homoclinic orbits or heteroclinic cycles. It then follows from  
\eeref{D+2}, \eeref{L0} and the Poincar\'e-Bendixson Theorem (see, for example, \cite{Verhulst90}), that
 $\mathcal{S}_1^+|_{v}$ must leave $\Omega_0$ through $L_2$ (at finite $y$) or connect with $\vec{e_a}$, for some $a \in [0,u_c]$ (as $y \to \infty$). For each $v \geq 0$, this observation allows us to classify the behaviour of $\mathcal{S}_1^+|_{v}$, by introducing the following function $\overline{s}: \overline{\mathbb{R}}^+ \to \mathbb{R}$, such that,
 \begin{center}
\begin{minipage}{0.85\linewidth}
	$\overline{s}(v)  = $ 
	The distance, measured from the origin of the $(\alpha, \beta)$ plane, to the point of intersection of $\mathcal{S}_1^+|_v$ with $L_2$  (negative distance)
	or $L_3$ (positive distance).
\end{minipage}	
\end{center}

\noindent  We have immediately that
\begin{equation}   \eelab{min_distancefunction}
\overline{s}(0)=\beta_0(0)<0,\quad\text{and}\quad
\beta_0(0)<\overline{s}(v) \leq u_c ,
\end{equation}
for all $v > 0$.\
Moreover, since $\vec{Q}(\alpha,\beta)$ depends continuously on $(\alpha,\beta,v) \in \overline{D}_+ \times \overline{\mathbb{R}}^+ \setminus  \{ (\beta,u_c) :  \beta \leq 0 \} \times \overline{\mathbb{R}}^+ $, the refraction vector \eeref{refraction vector} for phase paths crossing the line $\alpha = u_c$ in $D_+$ is independent of $(\beta,v) \in \mathbb{R}^- \times \overline{\mathbb{R}}^+$, and $\overline{\Omega}_0$ is compact, we may conclude that $\overline{s} \in C  (\overline{\mathbb{R}}^+)$. 
In addition, from the rotational properties of the vector field $\vec{Q}(\alpha , \beta)$ in $R^+$ with increasing $v \geq 0$, we deduce that
$\overline{s}(v_2) > \overline{s}(v_1)$ $\forall  v_2 > v_1 \geq 0$.
Therefore, $\overline{s}: \overline{\mathbb{R}}^+ \to \mathbb{R}$ is a continuous and strictly monotone increasing function. Next, take
\begin{equation}
v > v_c(u_c) = \left( \frac{1}{u_c} \sup_{\gamma \in (u_c,1]} f_c(\gamma) \right)^{\frac{1}{2}}   .
\end{equation}
Then, with $\beta_c = - v u_c$, we have
\begin{linenomath}
	\begin{equation}
\vec{Q}(\alpha , \beta_c) \cdot (0,1)  = v^2 u_c - f_c(\alpha) 
 >  \sup_{\gamma \in (u_c,1]} f_c(\gamma) -   f_c(\alpha)  \geq 0 , \eelab{phase_plane_direction}
\end{equation}
\end{linenomath}
for all $\alpha \in (u_c,1]$, and recall that $\mathcal{S}_0|_v$ is given by $\beta = - v \alpha$ for $\alpha \in [0,u_c]$. It then follows, from \eeref{phase_plane_direction}, that
\begin{equation} \eelab{s_u_geq_u_c}
\overline{s}(v) > 0 \quad \forall  v > v_c(u_c).
\end{equation}
We now observe that, at any $v \geq 0$, the dynamical system \eeref{phaseplane} has a heteroclinic connection between $\vec{e_1}$ and $\vec{e_0}$, in $R^+$ (which is unique, and is, in fact, contained in $\Omega_0 \subset R^+$) if and only if $\overline{s}(v)=0$. It follows that since $\overline{s} :\overline{\mathbb{R}}^+ \to \mathbb{R}$ is a continuous and strictly monotone increasing function, which satisfies \eeref{min_distancefunction} and \eeref{s_u_geq_u_c}, then, for each $u_c \in (0,1)$, there exists a unique $v^*(u_c) > 0$ such that
\begin{equation} \eelab{speed_zero_distance}
\overline{s}(v^*(u_c)) = 0 ,
\end{equation}
whilst,
\beq
\overline{s}(v) < 0 \quad \forall v \in [0, v^*(u_c)), \quad\text{and}\quad
\overline{s}(v) > 0 \quad \forall  v \in (v^*(u_c), \infty) .
  \eeq
We conclude that, for each $u_c \in (0,1)$, QIVP has a PTW solution if and only if $v=v^*(u_c)(>0)$ which we write as $u=U_T(y)$, $y \in \mathbb{R}$. Moreover, this PTW solution is unique. In addition, since the associated heteroclinic connection between $\vec{e_1}$ and $\vec{e_0}$ is contained in $\Omega_0$, then we conclude that $U_T: \mathbb{R} \to \mathbb{R}$ satisfies:
\begin{subequations}
\begin{equation} \eelab{u_domain}
0 < U_T(y) < 1, \quad U_T'(y) < 0 \qquad \forall  y \in \mathbb{R},
\end{equation}
with $U_T(0) = u_c$, and
 \begin{align}
 & U_T''(0^+)  -U_T''(0^-)  = - f_c^+ , \eelab{u_second_deriv_discontinuity} \\
 & U_T(y) = u_c e^{- v^*(u_c)y} \quad \forall  y \mbox{ } \in \overline{\mathbb{R}}^+, \eelab{u_cut-off} \\
 & U_T(y)  \sim 1 - A_{-\infty} e^{\lambda_+ ( v^*(u_c))y} \quad \mbox{as} \quad y \to - \infty, \eelab{u_small_y}
 \end{align}
\end{subequations}
for some constant $A_{-\infty}> 0$ (depending upon $u_c \in (0,1)$), and with the eigenvalue $\lambda_+(v)$ given in \eeref{evals_e1}.

We next consider $u_c \in (0,1)$ as a parameter, regarding $v^*$ as a function of $u_c$, with $v^*:(0,1) \to \mathbb{R}^+$ such that $v^*=v^*(u_c)$, and associated PTW solution $u=U_T(y,u_c)$ for $(y, u_c) \in \mathbb{R} \times (0,1)$. We recall that the vector field $\vec{Q}(\alpha, \beta)$ is continuously differentiable on
 {$(\alpha, \beta, v) \in ([0,u_c)\cup (u_c,1]) \times \mathbb{R} \times  \overline{\mathbb{R}}^+$},
whilst the refraction vector \eeref{refraction vector} depends on $u_c \in (0,1)$ and is continuous. It follows
that on fixing $u_c^0 \in (0,1)$, and taking $\varepsilon > 0$, then with $u_c=u_c^0$ and $v=v^*(u_c^0) - \varepsilon$,
we have that  $\overline{s}(v^*(u_c^0)-\varepsilon)|_{u_c=u_c^0} < 0$,
where we have used equation \eeref{speed_zero_distance}.
Hence, there exists $\delta_{\varepsilon}^- > 0$, which depends on $\varepsilon>0$, such that for all $u_c \in (u_c^0 - \delta_{\varepsilon}^- , u_c^0 + \delta_{\varepsilon}^-)= I^-_{\varepsilon}$,  {we have
$\overline{s}(v^*(u_c^0)-\varepsilon)|_{u_c \in  I^-_{\varepsilon}} < 0$}.
It follows that $v^*(u_c) > v^*(u_c^0) - \varepsilon$ for all $u_c \in I^-_{\varepsilon}$. Similarly, we establish that there exists $\delta_{\varepsilon}^+ > 0$, which depends on $\varepsilon >0$, such that for all $u_c \in (u_c^0 - \delta_{\varepsilon}^+ , u_c^0 + \delta_{\varepsilon}^+)= I^+_{\varepsilon}$,  {we have
$\overline{s}(v^*(u_c^0)+\varepsilon) |_{u_c \in I^+_{\varepsilon}} > 0$}.
It follows that $v^*(u_c) < v^*(u_c^0) + \varepsilon$ for all $u_c \in I^+_{\varepsilon}$. We now set $\delta_{\varepsilon} = \min (\delta_{\varepsilon}^-,\delta_{\varepsilon}^+)$. Thus, for all $u_c \in (u_c^0 - \delta_{\varepsilon} , u_c^0 + \delta_{\varepsilon})= I_{\varepsilon}$,  {we have 
$|  v^*(u_c) - v^*(u_c^0)  | < \varepsilon$}.
We conclude that $v^*:(0,1) \to \mathbb{R}$ is continuous. In addition, we recall that
\begin{equation} \eelab{v*pos}
v^*(u_c)>0 \quad \forall  u_c \in (0,1).
\end{equation}
Next, let $u_c^0 \in (0,1)$ and consider $\mathcal{S}_1^+ |_{(u_c^0,v^*(u_c^0))}$. It follows from the refraction vector \eeref{refraction vector} that there exists $\delta > 0$, such that on fixing $v=v^*(u_c^0)$, then for any $u_c \in (u_c^0, u_c^0 + \delta)=P_{\delta}$, the intersection point of $\mathcal{S}_1^+ |_{(u_c,v^*(u_c^0))}$ with the line $\alpha = u_c$ lies above the intersection point of the line $\beta = - v^*(u_c^0) \alpha$ with the line $\alpha = u_c$. Consequently, $\overline{s}(v^*(u_c^0)) |_{u_c \in P_{\delta} } > 0$, from which we conclude that $v^*(u_c) < v^*(u_c^0)$ for all $u_c \in P_{\delta}$.
Thus, $v^*:(0,1) \to \mathbb{R}$ is locally decreasing, and continuous, and so $v^*:(0,1) \to \mathbb{R}$ is strictly monotone decreasing. It then also follows from \eeref{v*pos} that $v^*(u_c)$ has a finite non-negative limit as $u_c \to 1^-$. Hence, $v^*(u_c) \to v_1^*$ as $u_c \to 1^-$, for some $v_1^* \geq 0$. When $(1-u_c)$ is sufficiently small, the linearisation theorem (see, for example, \cite{Verhulst90}) guarantees that $\mathcal{S}_1^+$ can be approximated in the region $(\alpha , \beta) \in [u_c, 1] \times \mathbb{R}^-$ by its linearised form at the equilibrium point $\vec{e_1}$; it is then readily established that $v_1^* = 0$, and, moreover, that  {$v^*(u_c) \sim  | f_c'(1) |^{\frac{1}{2}} (1 - u_c)$
as $u_c \to 1^- $}.
We now investigate $v^*(u_c)$ as $u_c \to 0^+$. To begin with we consider the dynamical system \eeref{phaseplane} when $u_c=0$.
In this case,
  the dynamical system \eeref{phaseplane} has a (unique) heteroclinic connection which connects $\vec{e_1}$, as $y \to - \infty$, to $\vec{e_0}$, as $y \to  \infty$, if and only if $v \in [2 , \infty)$, see for example
  \cite{Kolmogorov_etal1937,AronsonWeinberger1975, Larson1978, Needham1992}.
Moreover, $\overline{s}(v) |_{u_c=0} < 0$ for all $v \in [0, 2)$.
From \eeref{refraction vector} and \eeref{refraction_vector_limit}, it follows that $\mathcal{S}_1^+$ depends continuously on $u_c\geq 0$.
Thus, for $\varepsilon>0$, there exists $\sigma_{\varepsilon}>0$ such that for $u_c \in (0,\sigma_{\varepsilon})$, then $\overline{s}(2 - \varepsilon) |_{u_c} < 0$. Therefore, from \eeref{speed_zero_distance}, we deduce that $v^*(u_c) > 2 - \varepsilon$ for all $u_c \in (0,\sigma_{\varepsilon})$. However, it also follows from \eeref{refraction vector} and \eeref{refraction_vector_limit} that $ \overline{s}(2) |_{u_c} > 0$ for all $u_c \in (0,1)$.  Thus, $v^*(u_c) < 2$ for all $u_c \in (0,1)$. We conclude that,
 {$2 - \varepsilon < v^*(u_c) < 2$ $\forall \; u_c \in (0,\sigma_{\varepsilon})$.}
Since 
 {this} 
holds for all $\varepsilon >0$, we conclude immediately that $v^*(u_c)$ has limit $2$ as $u_c \to 0^+$.
We conclude that $v^*:(0,1) \to \mathbb{R}$ is continuous and monotone decreasing, with
\begin{equation}  \eelab{v*bound2}
\lim_{u_c \to 1^-} v^*(u_c) = 0, \qquad \lim_{u_c \to 0^+} v^*(u_c) = 2.
\end{equation}
This completes the proof of
Theorem \ref{theorem}.
In the next two sections we consider the structure of the PTW solutions in the limits $u_c \to 0^+$ and $u_c \to 1^-$ respectively.

\section{Asymptotic Structure of the PTW Solution when $u_c \to 0^+$}

In this section we investigate the detailed asymptotic form of $v^*(u_c)$ as $u_c \to 0^+$, in the small cut-off limit, via the method of matched asymptotic expansions. To that end, we write $u_c =\varepsilon$ with $ 0 < \varepsilon \ll 1$. It then follows from Theorem \ref{theorem} that we may write,
\begin{equation} \eelab{speed_smalluc1}
v^*(\varepsilon) = 2 - \bar{v}(\varepsilon),
\end{equation}
where now,
\begin{equation} \eelab{speed_smalluc2}
\bar{v}(\varepsilon) > 0 \quad \forall \; \varepsilon \in (0,1),
\end{equation}
and
\begin{equation} \eelab{speed_smalluc3}
\bar{v}(\varepsilon) = o(1) \quad \mbox{as} \quad \varepsilon \to 0^+.
\end{equation}
With $U_T: \mathbb{R} \to \mathbb{R}$ being the associated PTW solution, then from \eeref{BVP_cutoff},
\begin{linenomath}
\begin{subequations} \eelab{BVPsmalluc}
\begin{align}
& U_{Tyy} + (2 - \bar{v}(\varepsilon)) U_{Ty} + f(U_T) =0, \quad y<0, \\
& U_T(y) > \varepsilon \quad \forall \; y<0, \\
& U_T(0) = \varepsilon, \\
& U_{Ty}(0) = - (2- \bar{v}(\varepsilon) ) \varepsilon, \\
& U_T(y) \to 1 \quad \mbox{as} \quad y \to - \infty.
\end{align}
\end{subequations}
\end{linenomath}
It is convenient, in what follows, to make a shift of origin by introducing the coordinate $\bar{y}$ via
\begin{linenomath}
\begin{equation*}
\bar{y} = \bar{y}_c (\varepsilon) + y,
\end{equation*}
\end{linenomath}
where $\bar{y}_c (\varepsilon)$ is chosen so that \eeref{BVPsmalluc} becomes,
\begin{linenomath}
\begin{subequations} \eelab{BVP2smalluc}
\begin{align}
& U_{T\bar{y}\bar{y}} + (2 - \bar{v}(\varepsilon)) U_{T\bar{y}} + f(U_T(\bar{y})) =0, \quad \bar{y}<\bar{y}_c (\varepsilon), \\
& U_T(\bar{y}) > \varepsilon \quad \forall \; \bar{y}<\bar{y}_c (\varepsilon), \eelab{BVP2smalluca} \\
& U_T(\bar{y}_c (\varepsilon)) = \varepsilon, \eelab{BVP2smallucb} \\
& U_{T\bar{y}}(\bar{y}_c (\varepsilon)) = - (2- \bar{v}(\varepsilon) ) \varepsilon, \eelab{BVP2smallucc} \\
& U_T(\bar{y}) \to 1 \quad \mbox{as} \quad \bar{y} \to - \infty,
\end{align}
\end{subequations}
\end{linenomath}
with now the shift of origin fixing
\begin{equation} \eelab{fix_trans}
U_T(0) = \frac{1}{2}.
\end{equation}
It follows from \eeref{BVP2smalluc} and \eeref{fix_trans}
 that
\begin{equation} \eelab{y_smalluc}
\bar{y}_c (\varepsilon) \to + \infty \quad \mbox{as} \quad \varepsilon \to 0^+.
\end{equation}
Our objective is now to examine the boundary value problem \eeref{BVP2smalluc} and \eeref{fix_trans} as $\varepsilon \to 0^+$, and, in particular, to determine the asymptotic structure of $\bar{v}(\varepsilon)$ as $\varepsilon \to 0^+$. Anticipating the requirement of outer regions, we begin in an inner region when $\bar{y}=O(1)$ and $U_T=O(1)$ as $\varepsilon \to 0^+$, and we label this as region $\mathbf{I}$.
 In region $\mathbf{I}$ we thus expand as
\begin{equation} \eelab{exp_reg_I}
U_T(\bar{y} ; \varepsilon) = U_m(\bar{y}) + O ( \bar{v}(\varepsilon) ) \quad \mbox{as} \quad \varepsilon \to 0^+,
\end{equation}
with $\bar{y}=O(1)$. On substitution from \eeref{exp_reg_I} into \eeref{BVP2smalluc} and \eeref{fix_trans}, and using \eeref{y_smalluc}, we obtain the leading order problem as
\begin{linenomath}
\begin{subequations} \eelab{BVP_um_smalluc}
\begin{align}
& U_{m\bar{y}\bar{y}} + 2 U_{m\bar{y}} + f(U_m) =0, \quad - \infty < \bar{y} < \infty, \\
& U_m(\bar{y}) > 0, \quad - \infty < \bar{y} < \infty, \\
& U_m(\bar{y}) \to \begin{cases}
			1, & \text{as $\bar{y} \to - \infty$}\\
            0, & \text{as $\bar{y} \to \infty$} \end{cases} \\
& U_m(0) =\frac{1}{2}.
\end{align}
\end{subequations}
\end{linenomath}
The leading order problem is immediately recognised as the boundary value problem
\eeref{BVP_cutoff}
for the permanent form travelling wave solution to the corresponding KPP problem without cut-off $(\varepsilon =0)$.
Let $U_m:\mathbb{R} \to \mathbb{R}$ be the unique solution to \eeref{BVP_um_smalluc}. For use in what follows, we recall
\eeref{PTW_KPP} with higher order corrections given by
\begin{equation} \eelab{KPPextremes}
U_m(\bar{y}) =
\begin{cases}
			\left( A_{\infty} \bar{y} + B_{\infty} \right) e^{ - \bar{y}}+
				 O(\bar{y}^2e^{-2\bar{y}}), & \text{as $\bar{y} \to \infty$}\\
            1 - A_{- \infty} e^{\gamma \bar{y}}+O(e^{2\gamma \bar{y}}), & \text{as $\bar{y} \to -\infty$}
\end{cases}
\end{equation}
where
$\gamma = -1+\sqrt{1 + |f'(1)|}\quad (> 0)$.
On proceeding to $O(\bar{v}(\varepsilon))$ in region $\mathbf{I}$ 
we observe that the inner region expansion \eeref{exp_reg_I} becomes non-uniform when $| \bar{y} | \gg 1$, and in particular when $(-\bar{y})=O(\bar{v}(\varepsilon)^{- \frac{1}{2}} )$ and $\bar{y}=O(\bar{v}(\varepsilon)^{- \frac{1}{2}} )$. Therefore, to complete the asymptotic structure of the solution to \eeref{BVP2smalluc} as $\varepsilon \to 0^+$, we must introduce two outer regions, namely region $\mathbf{II^+}$ when $\bar{y}=O(\bar{v}(\varepsilon)^{- \frac{1}{2}} )$ and region $\mathbf{II^-}$ when $(-\bar{y})=O(\bar{v}(\varepsilon)^{- \frac{1}{2}} )$.
 {In this context, for any variable $\lambda$, we will henceforth write $\lambda=O(1)>0$ as $\lambda=O(1)^+$, and correspondingly, $\lambda=O(1)<0$ as $\lambda=O(1)^-$}. 
 We begin in region $\mathbf{II^-}$. 
 To formalize region $\mathbf{II^-}$, we introduce the scaled variable,
\begin{equation} \eelab{var_regionII}
\hat{y} = \bar{v}(\varepsilon)^{\frac{1}{2}} \bar{y},
\end{equation}
so that $\hat{y} = O(1)^-$ 
in region $\mathbf{II^-}$ as $\varepsilon \to 0^+$. It then follows from \eeref{exp_reg_I} and \eeref{KPPextremes} that
\begin{equation}
U_T( \hat{y} ; \varepsilon) = 1 - O \left( e^{- \bar{v}(\varepsilon)^{- \frac{1}{2}}}  \right),
\end{equation}
as $\varepsilon  \to 0^+$ in region $\mathbf{II^-}$. It is then straightforward to develop an exponential expansion in region $\mathbf{II^-}$, which, after matching (following the Van Dyke matching principle, \cite{VanDyke1975}) with region $\mathbf{I}$, via \eeref{exp_reg_I} and \eeref{KPPextremes}, gives the outer expansion in region $\mathbf{II^-}$ as,
\begin{linenomath}
\begin{align}
 U_T( \hat{y} ; \varepsilon) = 1& - A_{- \infty}
 \exp\left[\gamma \bar{v}(\varepsilon)^{- \frac{1}{2}} \left(  1 + O(\bar{v}(\varepsilon) ) \right) \hat{y}    \right]  \nonumber\\
 &+O \left( \exp\left[ 2 \gamma \bar{v}(\varepsilon)^{- \frac{1}{2}} \left(  1 + O(\bar{v}(\varepsilon) ) \right) \hat{y}    \right]  \right),\eelab{outer_IIL}
\end{align}
\end{linenomath}
as $\varepsilon \to 0^+$ with $\hat{y} = O(1)^-$.
Thus, the solution in region $\mathbf{II^-}$ is at this order unaffected  by the cut-off.
We now proceed to region $\mathbf{II^+}$, where $\hat{y} = O(1)^+$ as $\varepsilon \to 0^+$. It is within this region that the conditions at $\bar{y} = \bar{y}_c (\varepsilon)$ must be satisfied, which then requires $\bar{y}_c (\varepsilon) = O ( \bar{v}(\varepsilon)^{- \frac{1}{2}} )$ as $\varepsilon \to 0^+$, which is consistent with \eeref{y_smalluc}. Thus, we write
\begin{equation} \eelab{outerIIR}
\bar{y}_c (\varepsilon) = \bar{v}(\varepsilon)^{- \frac{1}{2}} \hat{y}_c (\varepsilon),
\end{equation}
so that now,
\begin{equation} \eelab{outerIIRb}
\hat{y}_c (\varepsilon) = O(1)^+ \quad \mbox{as} \quad \varepsilon \to 0^+.
\end{equation}
In region $\mathbf{II^+}$ it follows from \eeref{exp_reg_I} and \eeref{KPPextremes} that
\begin{linenomath}
\begin{equation*}
U_T( \hat{y} ; \varepsilon) =O \left( \bar{v}(\varepsilon)^{- \frac{1}{2}} e^{- \bar{v}(\varepsilon)^{- \frac{1}{2}}}  \right),
\end{equation*}
\end{linenomath}
as $\varepsilon \to 0^+$. Again, it is then straightforward to develop an exponential expansion in region $\mathbf{II^+}$, which, after matching with region $\mathbf{I}$, via \eeref{exp_reg_I} and \eeref{KPPextremes}, gives the outer expansion in region $\mathbf{II^+}$ as,
\begin{linenomath}
\begin{align}
U_T( \hat{y} ; \varepsilon) =
&\left( A_{\infty} \bar{v}(\varepsilon)^{- \frac{1}{2}} \sin \big( \hat{y} ( 1 + O (\bar{v}(\varepsilon))) \big)  + B_{\infty}  \cos \big( \hat{y} ( 1 + O (\bar{v}(\varepsilon))) \big)              \right)   \nonumber\\
& \times \exp\left[- \bar{v}(\varepsilon)^{- \frac{1}{2}} \left( 1 + O(\bar{v}(\varepsilon) )    \right) \hat{y}\right]   \nonumber\\
&+ O \left(   \exp\left[ - 2\bar{v}(\varepsilon)^{- \frac{1}{2}} \left( 1 + O(\bar{v}(\varepsilon) )    \right) \hat{y}\right]        \right)\eelab{outer_IIR}
\end{align}
\end{linenomath}
as $\varepsilon \to 0^+$ with $\hat{y} = O(1)^+$. It now remains to apply conditions \eeref{BVP2smalluca}, \eeref{BVP2smallucb} and \eeref{BVP2smallucc} to \eeref{outer_IIR}. In the outer region $\mathbf{II^+}$, these conditions become,
\begin{linenomath}
\begin{subequations} \eelab{conditionsII}
\begin{align}
& U_T(\hat{y}; \varepsilon) > \varepsilon \quad \forall\;   O \left( \bar{v}(\varepsilon)^{- \frac{1}{2}} \right)^+ < \hat{y}<\hat{y}_c (\varepsilon), \eelab{conditionsIIa} \\
& U_T(\hat{y}_c (\varepsilon); \varepsilon) = \varepsilon, \eelab{conditionsIIb} \\
& U_{T\hat{y}}(\hat{y}_c (\varepsilon); \varepsilon) = - \varepsilon \bar{v}(\varepsilon)^{- \frac{1}{2}} (2- \bar{v}(\varepsilon) ) . \eelab{conditionsIIc}
\end{align}
\end{subequations}
\end{linenomath}
We now turn to conditions \eeref{conditionsIIb} and \eeref{conditionsIIc}. It is convenient to first eliminate $\varepsilon$ explicitly between \eeref{conditionsIIb} and \eeref{conditionsIIc} to give,
\begin{equation} \eelab{eq:con4 on uT region II right chap 4}
U_{T \hat{y}} ( \hat{y}_c (\varepsilon); \varepsilon ) = - \bar{v}(\varepsilon)^{- \frac{1}{2}} (2- \bar{v}(\varepsilon) ) U_T(\hat{y}_c (\varepsilon); \varepsilon),
\end{equation}
which replaces \eeref{conditionsIIc}. On substitution from \eeref{outer_IIR} into \eeref{eq:con4 on uT region II right chap 4} and expanding, using \eeref{speed_smalluc2}, \eeref{speed_smalluc3} and \eeref{outerIIRb}, we obtain,
\begin{linenomath}
\begin{align}
    A_{\infty} & \sin \omega   =- \bar{v}(\varepsilon)^{\frac{1}{2}}   (A_{\infty} + B_{\infty} )
    \cos \omega,\qquad \omega= \hat{y}_c (\varepsilon) ( 1 + O (\bar{v}(\varepsilon)) )
    \eelab{eq:region II right problem}
\end{align}
\end{linenomath}
as $\varepsilon \to 0^+$. Following \eeref{outerIIRb} and \eeref{eq:region II right problem}, we now expand,
\begin{equation} \eelab{eq:region II right expansion chap 4}
\hat{y}_c ( \varepsilon) = \hat{y}_c^0 + \hat{y}_c^1 \bar{v}(\varepsilon)^{\frac{1}{2}} +  O(\bar{v}(\varepsilon) ),
\end{equation}
as $\varepsilon \to 0^+$, with the constants $\hat{y}_c^0\;(>0)$ and $\hat{y}_c^1$ to be determined. On substitution from \eeref{eq:region II right expansion chap 4} into \eeref{eq:region II right problem}, we obtain, at $O(1)$,
\begin{linenomath}
\begin{equation*}
A_{\infty} \sin \hat{y}_c^0 = 0.
\end{equation*}
\end{linenomath}
\begin{figure}
\centerline{\includegraphics[width=0.99\linewidth]{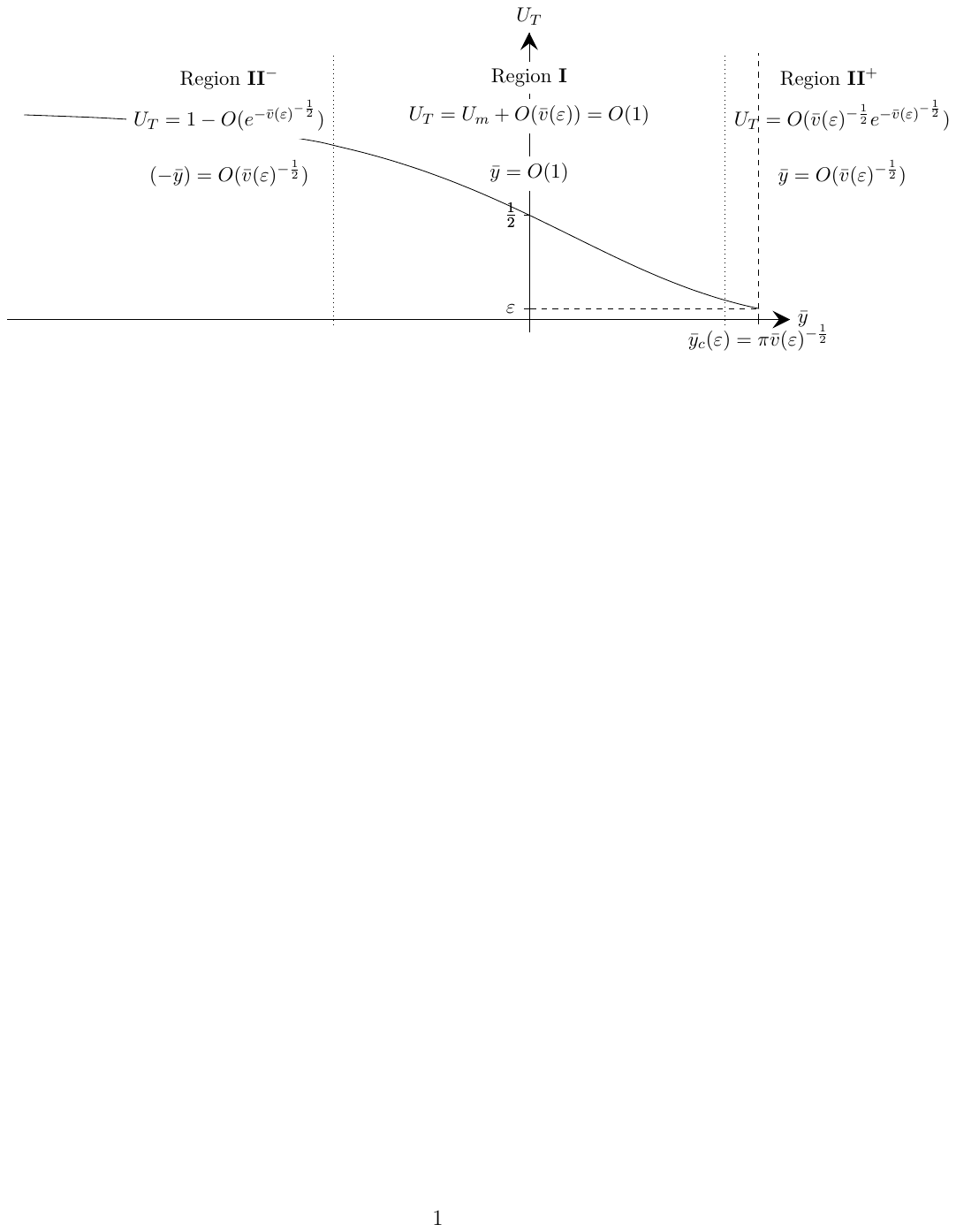}}
\caption{ {A schematic diagram of the asymptotic structure of  $U_T(\bar y;\varepsilon)$   as $u_c\to 0^+$}.
}
\label{fig:sketch}
\end{figure}
Since $A_{\infty} >0$, then we must have (recalling $\hat{y}_c^0>0$) $\hat{y}_c^0=k \pi$, for some $k \in \mathbb{N}$. However, condition \eeref{conditionsIIa}, with \eeref{outer_IIR}, then requires $k=1$, and so
\begin{equation} \eelab{eq:y hat c leading order chap 4}
\hat{y}_c^0 = \pi.
\end{equation}
Proceeding to $O(\bar{v}(\varepsilon)^{\frac{1}{2}} )$, we find that, on using \eeref{eq:y hat c leading order chap 4},
\begin{equation} \eelab{eq:y hat c second order chap 4}
\hat{y}_c^1= - \frac{(A_{\infty} + B_{\infty} )}{A_{\infty}} .
\end{equation}
Thus, via \eeref{eq:region II right expansion chap 4}, \eeref{eq:y hat c leading order chap 4} and \eeref{eq:y hat c second order chap 4} we have,
\begin{equation} \eelab{regionII_exp}
\hat{y}_c ( \varepsilon) = \pi -  \frac{(A_{\infty} + B_{\infty} )}{A_{\infty}} \bar{v}(\varepsilon)^{\frac{1}{2}} + 
O(\bar{v}(\varepsilon) ),
\end{equation}
as $\varepsilon \to 0^+$. It remains to apply condition \eeref{conditionsIIb}. On using \eeref{outer_IIR} and \eeref{regionII_exp}, condition \eeref{conditionsIIb} becomes
\begin{equation} \eelab{vbar}
\ln \varepsilon = - \frac{\pi}{\bar{v}(\varepsilon)^{\frac{1}{2}}} + \left(   \frac{(A_{\infty} + B_{\infty} )}{A_{\infty}} +  \ln A_{\infty} \right) + 
O(\bar{v}(\varepsilon)^{\frac{1}{2}}), 
\end{equation}
as $\varepsilon \to 0^+$. A re-arrangement of \eeref{vbar} then gives,
\begin{equation} \eelab{vbar2}
\bar{v}(\varepsilon)= \frac{\pi^2}{(\ln \varepsilon)^2} + \frac{ 2 \pi^2 \left(   (A_{\infty} + B_{\infty} )A_{\infty}^{-1} + \ln A_{\infty}       \right)}{ ( \ln \varepsilon)^3} + 
O \left(  \frac{1}{ ( \ln \varepsilon^4)}   \right), 
\end{equation}
as $\varepsilon \to 0^+$. It then follows from \eeref{regionII_exp} and \eeref{vbar2} that,
\begin{equation} \eelab{eq:y hat determined chap 4}
\hat{y}_c (\varepsilon) = \pi + \frac{(A_{\infty} + B_{\infty} ) \pi }{A_{\infty}} \frac{1}{ \ln \varepsilon} +   O \left(  \frac{1}{ ( \ln \varepsilon^2)}   \right),
\end{equation}
as $\varepsilon \to 0^+$. Finally, via \eeref{speed_smalluc1} and \eeref{vbar2}, we can construct $v^*(\varepsilon)$ as,
\begin{equation} \eelab{speed_small_cutoff}
v^*(\varepsilon)= 2 - \frac{\pi^2}{(\ln \varepsilon)^2} - \frac{ 2 \pi^2 \left(    (A_{\infty} + B_{\infty} )A_{\infty}^{-1} + \ln A_{\infty}       \right)}{ ( \ln \varepsilon)^3} + 
O \left(  \frac{1}{ ( \ln \varepsilon^4)}   \right), 
\end{equation}
as $\varepsilon \to 0^+$.
 {For completeness we give a schematic diagram of the asymptotic structure for
$U_T(\bar y;\varepsilon)$ in terms of the coordinate $\bar y$ as $\varepsilon\to 0^+$
in Figure \ffref{sketch}. 
Returning to  \eeref{speed_small_cutoff} 
}
we  observe that the   approximation   is decreasing in $\varepsilon$ as $\varepsilon \to 0^+$, and is in full accord with the rigorous results established in Theorem \ref{theorem}.
We see immediately that  
the approximation derived here,
agrees in the first two terms with  prediction \eeref{speed_cutoff_BD} that Brunet and Derrida \cite{BrunetDerrida1997} first obtained and its third term is consistent with the order of the error term in \eeref{speed_cut-off_Dumortier_etal}
that Dumortier, Popovic and Kaper  \cite{Dumortier_etal2007} derived.
   However, the method of matched asymptotic expansions 
    has enabled us to obtain the next correction term in \eeref{speed_small_cutoff}, and higher order terms could be obtained by systematically following this approach (of course it may also be possible  to obtain the third and higher order terms  via extending the approach of  \cite{Dumortier_etal2007}). 
 In fact,  
  we may continue the expansion in each region to next order, and after matching, we can readily obtain that 
 the higher-order correction to \eeref{speed_small_cutoff}  is given by 
 \begin{equation}\eelab{speed_small_cutoffb} 
 	\begin{split}
 v^*(\varepsilon)= 2 &- \frac{\pi^2}{(\ln  \varepsilon)^2} - \frac{ 2 \pi^2 \left(    (A_{\infty} + B_{\infty} )A_{\infty}^{-1} + \ln A_{\infty}       \right)}{ ( \ln  u_c)^3}\\ +
 &\frac{  3\pi^2 \left(\frac{1}{4}\pi^2-\left(    (A_{\infty} + B_{\infty} )A_{\infty}^{-1} + \ln A_{\infty}       \right)^2\right)}{ ( \ln  \varepsilon)^4}
+ O \left(  \frac{1}{ ( \ln  \varepsilon)^5}   \right),
 \end{split}
 \end{equation}
   as $\varepsilon\to 0^+$. For brevity, we do not provide a derivation to  \eeref{speed_small_cutoffb}. 
We now consider the asymptotic structure of the PTW solution to QIVP as $u_c \to 1^-$.

\section{Asymptotic Structure of the PTW Solution when $u_c \to 1^-$}

In this section we investigate the asymptotic form of $v^*(u_c)$ in the large cut-off limit $u_c \to 1^-$. To this end, we write $u_c =1 - \delta$ with $ 0 < \delta \ll 1$. Theorem \ref{theorem} guarantees the existence and uniqueness of a PTW solution, whose speed $v^*(\delta)=o(1)$ as $\delta \to 0^+$.
In this case, it is most convenient to consider the problem in the $(\alpha, \beta)$ phase plane corresponding to the phase path representing the PTW when $u_c = 1 - \delta$ and $v=v^*(\delta)$. Via \eeref{phaseplane}, \eeref{evals_e1}, \eeref{unstable_manifold} and \eeref{stable manifold of e0}, this is given by the phase path $\beta=\beta(\alpha; \delta)$, which satisfies the boundary value problem
\begin{linenomath}
\begin{subequations} \eelab{DS_largeuc}
\begin{align}
& \frac{d \beta}{d \alpha}  = -v^*(\delta) - \frac{f(\alpha)}{\beta}, \quad  \alpha \in (1 - \delta , 1),  \\
& \beta (\alpha; \delta)  \sim  - \lambda_+(v^*(\delta))(1  -  \alpha) \quad  \mbox{as} \quad \alpha \to 1^-, \eelab{evec_e1_largeuc} \\
& \beta (1 - \delta; \delta)  =  - v^*(\delta) (1 - \delta).
\eelab{bc_largeuc}
\end{align}
\end{subequations}
\end{linenomath}
We now examine the boundary value problem \eeref{DS_largeuc} as $\delta \to 0^+$.
Since $v^*(\delta)=o(1)$ as $\delta \to 0^+$, we expand $\lambda_+(v^*(\delta))$, via \eeref{evals_e1}, which determines that $\lambda_+(v^*(\delta))=O(1)$ as $\delta \to 0^+$. It follows from the boundary condition \eeref{evec_e1_largeuc}, that $\beta = O(\delta)$ as $\delta \to 0^+$. We therefore introduce the following re-scalings
\begin{equation} \eelab{re-scale_largeuc}
\beta = \delta Y, \quad \alpha = 1 - \delta X,
\end{equation}
with $Y,X=O(1)$ as $\delta \to 0^+$. The form of the boundary condition \eeref{bc_largeuc} then necessitates that $v^*(\delta) = O(\delta)$ as $\delta \to 0^+$. Thus, we write
\begin{equation} \eelab{speed_largeuc}
v^*(\delta) = \delta V (\delta) ,
\end{equation}
where $V(\delta)=O(1)$ as $\delta \to 0^+$. These re-scalings transform the boundary value problem \eeref{DS_largeuc} into,
\begin{linenomath}
\begin{subequations}\eelab{DS2_largeuc}
\begin{align}
& \frac{d Y}{d X} = \delta V(\delta) + \frac{f(1-\delta X)}{\delta Y}, \quad X \in (0,1),  \eelab{DS2_largeuca}  \\
& Y(X; \delta)  \sim - \lambda_+ ( \delta V(\delta)) X \quad \mbox{as} \quad X \to 0^+, \\
& Y(1 ; \delta)  = - V(\delta) (1 - \delta).
\end{align}
\end{subequations}
\end{linenomath}
We now expand $Y(X;\delta)$ and $V(\delta)$ according to,
\begin{linenomath}
\begin{subequations} \eelab{expansions_largeuc}
\begin{align}
& Y(X;\delta) = Y_0(X) + \delta Y_1(X) + o(\delta) ,   \quad X \in [0,1], \eelab{Yexpansion_largeuc}  \\
& V(\delta) = V_0 + \delta V_1 + o(\delta) , \eelab{Vexpansion_largeuc}
\end{align}
\end{subequations}
\end{linenomath}
as $\delta \to 0^+$.
Substituting the expansions from \eeref{expansions_largeuc} into the boundary value problem \eeref{DS2_largeuc} and expanding, at $O(1)$, we obtain the following boundary value problem for $Y_0(X)$, namely,
\begin{linenomath}
\begin{subequations}
\begin{align}
& \frac{d Y_0}{d X}  = -  f'(1) \frac{  X}{Y_0}, \quad X \in (0,1) , \eelab{Y0}   \\
& Y_0(X) \sim - | f'(1) |^{\frac{1}{2}} X \quad \mbox{as} \quad X \to 0^+, \eelab{bcY0} \\
& Y_0(1) = - V_0. \eelab{V0}
\end{align}
\end{subequations}
\end{linenomath}
The general solution to \eeref{Y0} is $Y_0^2(X) = c_1 -  f'(1) X^2$, for $X \in [0,1]$, where $c_1$ is an arbitrary constant of integration. Applying the boundary condition \eeref{bcY0} determines $c_1=0$. Therefore,
\begin{equation} \eelab{Y0soln}
Y_0(X) =  - | f'(1) |^{\frac{1}{2}} X, \quad  X \in [0,1].
\end{equation}
Application of the boundary condition \eeref{V0} then determines
\begin{equation} \eelab{V0b}
V_0 =   | f'(1) |^{\frac{1}{2}} .
\end{equation}
At $O(\delta)$, we obtain the following boundary value problem for $Y_1(X)$, namely,
\begin{linenomath}
\begin{subequations}
\begin{align}
& \frac{d Y_1}{d X} -   \frac{Y_1}{Y_0(X)^2} f'(1) X   = V_0 + \frac{1}{2} f''(1) \frac{ X^2  }{Y_0(X)}, \quad  X \in (0,1) ,
\eelab{Y1}   \\
& Y_1(X)  \sim \frac{1}{2}V_0 X  \quad \mbox{as} \quad X \to 0^+,
\eelab{Y1bc} \\
& Y_1(1)  = V_0 - V_1.  \eelab{V1a}
\end{align}
\end{subequations}
\end{linenomath}
On substituting $Y_0(X)$, given by \eeref{Y0soln}, into equation \eeref{Y1} and solving, we find that the general solution is
\begin{equation}
Y_1(X) = \frac{1}{2}V_0 X - \frac{1}{6} \frac{  f''(1) }{ | f'(1) |^{\frac{1}{2}}}X^2 + \frac{c_2}{X}, \quad  X \in (0,1],
\end{equation}
where $c_2$ is an arbitrary constant of integration. From the boundary condition \eeref{Y1bc}, $Y_1(X)$ remains bounded as $X \to 0^+$. Therefore, we require $c_2=0$. Thus, we obtain the solution for $Y_1(X)$ as
\begin{equation} \eelab{Y1soln}
Y_1(X)  = \frac{1}{6} | f'(1) |^{\frac{1}{2}}X \left( 3 -  \frac{ f''(1)  }{  | f'(1) |}  X  \right), \quad X \in [0,1].
 \end{equation}
Finally, an application of the boundary condition \eeref{V1} determines
\begin{equation}
V_1  = \frac{1}{6} | f'(1) |^{\frac{1}{2}}\left( 3 + \frac{ f''(1)   }{  | f'(1) |} \right) . \eelab{V1}
\end{equation}
On collecting expressions \eeref{Yexpansion_largeuc}, \eeref{Y0soln} and \eeref{Y1soln}, we have established that
\begin{linenomath}
\begin{align}
Y(X; \delta) = - | f'(1) |^{\frac{1}{2}} X & + \frac{1}{6} \delta  | f'(1) |^{\frac{1}{2}}X \left( 3 -  \frac{ f''(1) }{  | f'(1) |}   X  \right)  \nonumber \\
& + o(\delta) \quad \mbox{as} \quad \delta \to 0^+, \eelab{large_cut-offY}
\end{align}
\end{linenomath}
uniformly for $X \in [0,1]$. Similarly, on collecting expressions \eeref{Vexpansion_largeuc}, \eeref{V0b} and \eeref{V1}, we obtain,
\begin{equation}\eelab{Vdelta}
V(\delta) =  | f'(1) |^{\frac{1}{2}}  + \frac{1}{6} \delta | f'(1) |^{\frac{1}{2}}  \left( 3 + \frac{f''(1)}{ | f'(1) |} \right)+ o(\delta) \quad \mbox{as} \quad \delta \to 0^+.
\end{equation}
We use \eeref{re-scale_largeuc} to express the PTW solution to QIVP in terms of the cut-off $u_c$ as
\begin{linenomath}
\begin{align}
\beta(\alpha)
= - \frac{1}{2}| f'(1) |^{\frac{1}{2}}
(1+u_c) (1-\alpha) &
-\frac{1}{6}\frac{f''(1)}{| f'(1) |^{\frac{1}{2}}}(1-\alpha)^2
    \nonumber \\
& + o((1-u_c)^2) \quad \mbox{as} \quad u_c \to 1^-, \eelab{large_cut-offbeta}
\end{align}
\end{linenomath}
with $\alpha\in [u_c,1]$.
Its speed of propagation, via \eeref{speed_largeuc} and \eeref{Vdelta}, is given by
\begin{linenomath}
\begin{align}
v^*(u_c)  =(1 - u_c) | f'(1) |^{\frac{1}{2}} & + \frac{1}{6} (1 - u_c)^2 | f'(1) |^{\frac{1}{2}}\left( 3+ \frac{f''(1)}{| f'(1) |} \right) \nonumber \\
& + o((1 - u_c)^2)  \quad \mbox{as } u_c \to 1^-. \eelab{speed_large_cutoff}
\end{align}
\end{linenomath}
In the next section we consider the specific case of a cut-off Fisher reaction,  determining  $U_T:\mathbb{R}\to \mathbb{R}$ and $v^*:(0,1) \to \mathbb{R}$ via numerical integration.

\section{Numerical Example} \label{sec:Numerical example}

We now focus on the particular case of the cut-off Fisher reaction function, namely, \eeref{BDreaction} with \eeref{Fisher}.
We obtain numerical approximations of the speed $v^*(u_c)$  and PTW solutions $U_T:\mathbb{R}\to\mathbb{R}$
 for a range of values of cut-off $u_c$.
This is achieved
by solving \eeref{NBVP_KPPbintro}
numerically over an interval 
  {$y\in[0,M]$}
for $M\in\mathbb{R}^+$  using the Matlab initial value solver \texttt{ode45} where $U_T(0)=u_c$.
 As `initial' condition we
employ  \eeref{unstable_manifold} to approximate the unstable manifold
 near the unstable fixed point $(U_T,U_T')=(1,0)$, 
 taking $U_T(0)=1-\epsilon$ and $U_T'(0)=-\lambda_+(v)\epsilon$
 where $\epsilon=10^{-12}$ and prescribe an absolute and relative tolerance of $10^{-13}$.
 The  value of $v$ in the second initial condition is not known {\it a priori}.
 We therefore build the initial value solver into a shooting type algorithm, for which we guess the value of $v$, integrate \eeref{NBVP_KPPbintro}  
 to obtain $U_T'(0)$, and then compare $U_T'(0)$ to the target value $-v u_c$.
 The value of $v$ is then modified using the bisection method and this integration procedure is iterated until the absolute error  satisfies
  {$|U_T'(0)+vu_c|<10^{-13}$}.
We  start with a value of $u_c$ close to $1$ and take  $v=2$ as initial guess to begin the iteration which leads to  the solution  with speed $v^*(u_c)$.
We then   iterate over decreasing values of $u_c$ using the previously determined value of $v^*(u_c)$   as an initial guess to find the next solution.

\begin{figure}
 		     \centering
 		  	     \includegraphics[width=0.6\textwidth]{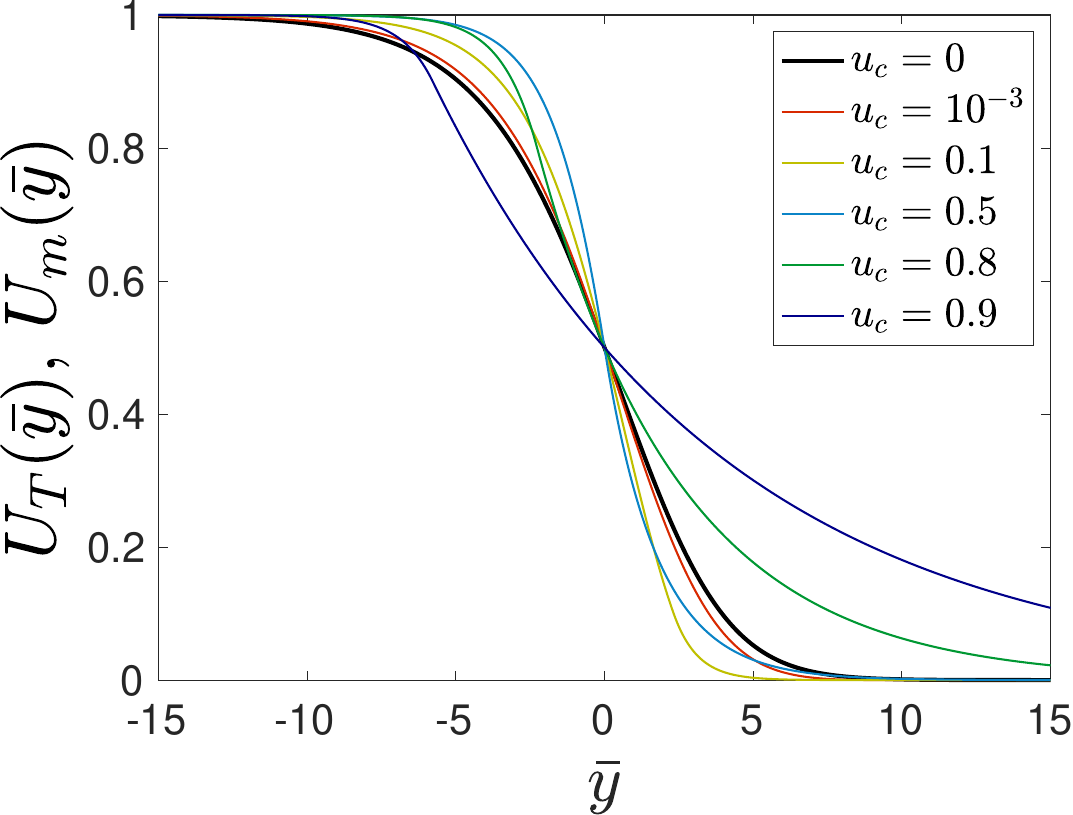}
\caption{
Permanent form travelling wave solutions
  $U_T(\bar{y})$ computed numerically as a function of $\bar{y}$ for  selected cut-off values $u_c$. These are compared against the  travelling wave solution $U_m(\bar{y})$ obtained in the absence of a cut-off (corresponding to $u_c=0$).
}
 \fflab{PTW}
 \end{figure}

It is also useful to obtain a numerical approximation of the permanent traveling wave solution $U_m$ for the
Fisher reaction function \eeref{Fisher}  in the absence of a cut-off.
This is readily achieved by
solving \eeref{NBVP_KPPintro}
numerically over an interval  {$y\in[0,N]$ for $N\in\mathbb{R}^+$}  once more using the Matlab initial value solver \texttt{ode45}.  
 As `initial condition'
we employ \eeref{PTW_KPPb} to approximate the unstable manifold
 near the unstable fixed point $(U_m,U_m')=(1,0)$,  
 taking $v=2$,  $U_m(0)=1-\epsilon$ and $U_m'(0)=(\sqrt{2}-1)\epsilon$
 where $\epsilon=10^{-12}$ and prescribe an absolute and relative tolerance of $10^{-13}$. We then determine that value of $y$ for which $U_m$ is equal to $1/2$ and then perform a coordinate shift to the origin.

 Figure \ffref{PTW} contrasts the  behaviour of  $U_T$    against $U_m$. %
A direct comparison between $U_T$ and $U_m$ is achieved when $U_T$ and $U_m$ are expressed in terms of $\bar{y}$ so that
  $U_T(0)= U_m(0)\approx 0.5$.
For small values of $u_c$, $U_T(\bar{y})$ is close to $U_m(\bar{y})$ in agreement with the asymptotic theory of section 4.  As $u_c$ increases, we observe a strong departure  {of $U_T(\bar{y})$} from $U_m(\bar{y})$ with the slope  of $U_T(\bar{y})$  at $\bar{y}=0$ becoming increasingly steep  until it reaches a maximum  at $u_c=0.5$. 
Beyond this value, the slope satisfies  {$U_T'(0)=-v^*(u_c)/2$}   
and therefore
  becomes increasingly gentle with $u_c$ (since $v^*(u_c)$ decreases with $u_c$).  
When $u_c$ approaches $1$,
$U_T(\bar{y})$ is in full agreement with
 the asymptotic prediction, derived from \eeref{large_cut-offbeta}  {(not shown)}.

   \begin{figure}[t]
     \centering
     \begin{minipage}[b]{0.48\textwidth}
  	     \includegraphics[width=\textwidth]{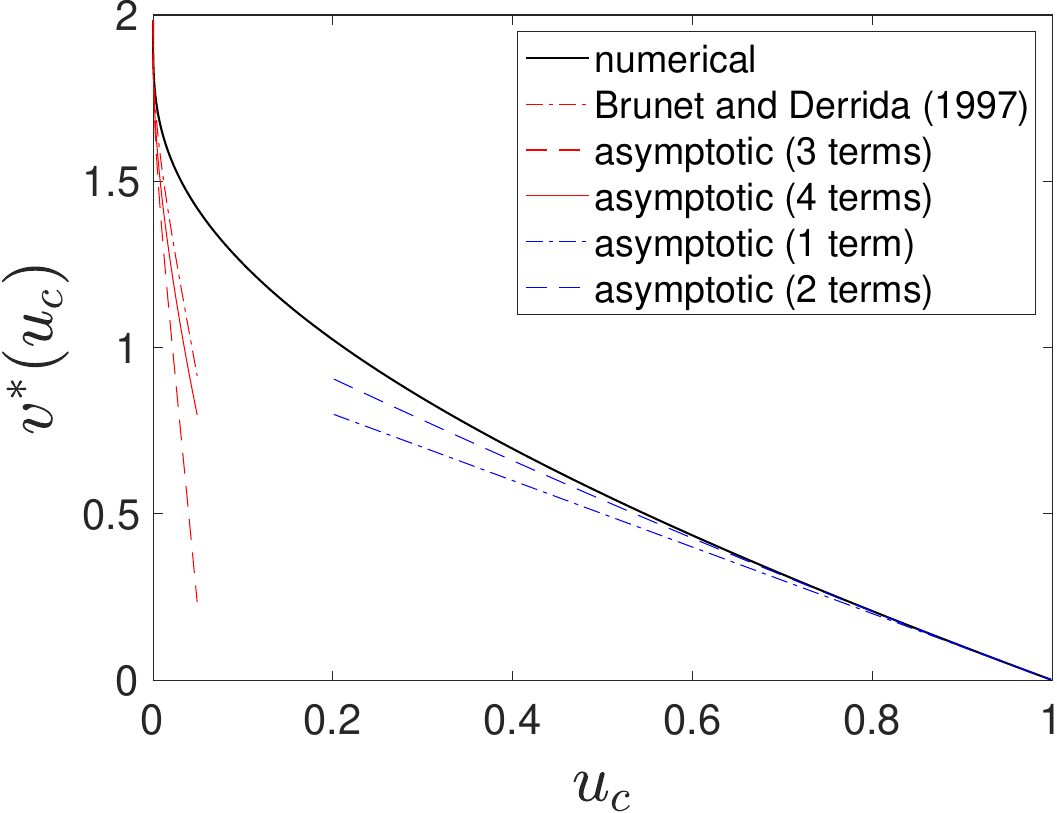}\\
  	     \centering  (a)
     \end{minipage}
     \hfill
     \begin{minipage}[b]{0.49\textwidth}
  	   \includegraphics[width=\textwidth]{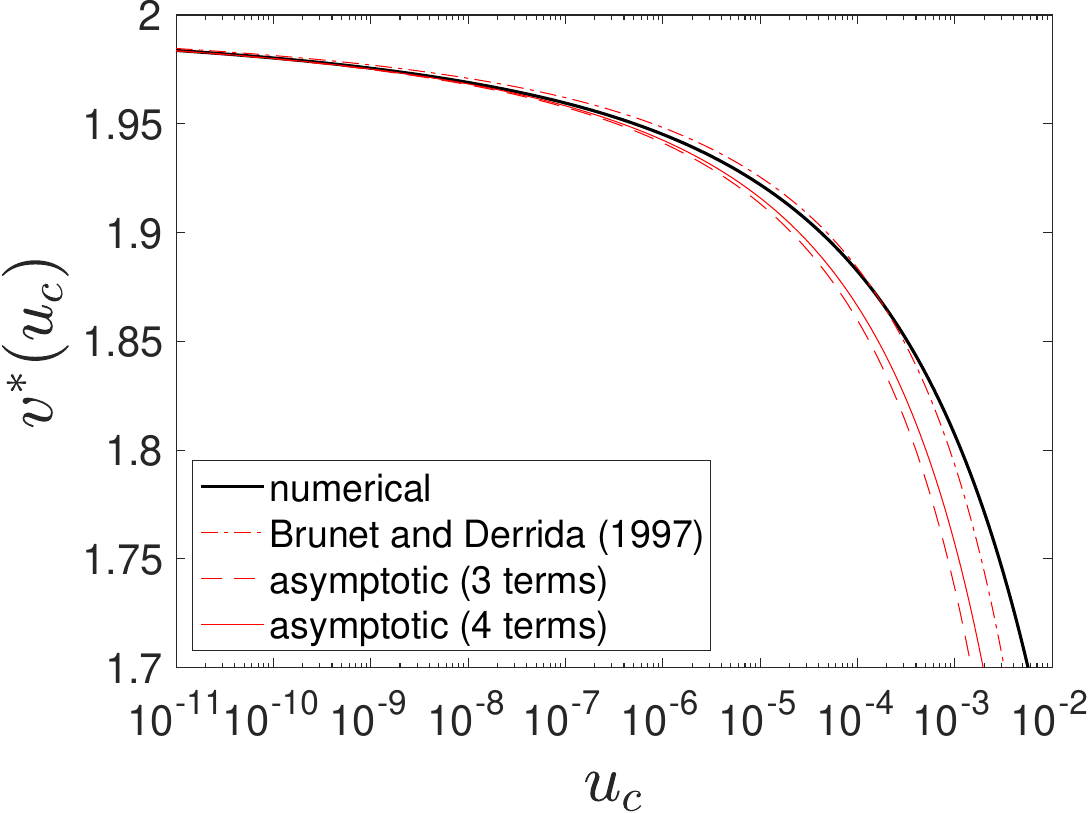}\\
    \centering  (b)
     \end{minipage}
   \caption{
     (a) Propagation speed $v^*(u_c)$ 
      computed numerically as a function of the cut-off value $u_c$ for  the particular case of the cut-off Fisher reaction function \eeref{BDreaction} with \eeref{Fisher}. Comparison against the   asymptotic expansions for $v^*(u_c)$ 
	 obtained  
  as  $u_c\to 1^-$ 
  based on retaining  one and two terms in \eeref{speed_large_cutoff} 
	  and    
  as $u_c\to 0^+$ based on retaining two, three and four terms in 
   \eeref{speed_small_cutoffb}. 
   The two-term asymptotic   expansion of $v^*(u_c)$ as $u_c\to 0^+$ corresponds to   prediction  \eeref{speed_cutoff_BD}  first  obtained in \cite{BrunetDerrida1997}.  
 (b) Same as (a) but focusing on smaller values of
 	    $u_c$.
  \fflab{speed}
  }
   \end{figure}

 Figure \ffref{speed}  examines the behaviour of the speed $v^*(u_c)$ and compares it to the various  asymptotic expansions  obtained as   $u_c\to 0^+$ and   $u_c\to 1^-$ %
  (as derived from \eeref{speed_small_cutoffb} and \eeref{speed_large_cutoff},  respectively).
  The  asymptotic  expansion  \eeref{speed_small_cutoffb} obtained for  $u_c\to 0^+$ relies on
   on  the global constants $A_\infty$ and $B_\infty$ 
    associated with the leading edge behaviour of
   $U_m(\bar{y})$
   (see \eeref{KPPextremes}).
   We determine the  values of $A_\infty$ and $B_\infty$
   by performing a least-squares polynomial
   fit to the computed $U_m(\bar{y}) e^{\bar{y}}$ for $\bar{y}\gtrsim 10$ from where we obtain   to a very good approximation the linear polynomial fit with
   \beq
   A_\infty\approx 3.5\quad\text{and}\quad B_\infty\approx-11.2.
   \eeq
  Figure \ffref{speed}(a)   demonstrates that the two-term asymptotic expansion  of
   $v^*(u_c)$ as $u_c\to 1^-$  accurately captures  the speed $v^*(u_c)$ for a wide range of values   given by $0.4\lesssim u_c<1$
  (when $u_c=0.4$, $\delta=1-u_c=0.6$ associated with expansions \eeref{expansions_largeuc}
  is no longer small). 
  Figure \ffref{speed}(b) focusses on the behaviour of the speed  obtained for smaller values of $u_c$.
   It shows that  the curve representing  
  the two-term expansion  based on retaining two terms in  \eeref{speed_small_cutoffb}  and corresponding to the  asymptotics that \cite{BrunetDerrida1997}  first obtained  
    crosses the numerically computed curve
   representing $v^*(u_c)$ at  $u_c\approx 10^{-4}$ 
    and has a monotonic approach to this curve for smaller values of $u_c$. We therefore anticipate that this two-term   expansion
    only becomes genuinely asymptotic for values of $u_c\ll 10^{-4}$. 
  This implies that any comparison of the   three-term and four-term 
  expansions  based on retaining three and four terms in \eeref{speed_small_cutoffb}, respectively, 
    should only be considered for  $u_c\ll 10^{-4}$. 
    With this in mind,  
    the logarithmic corrections included in 
   our three-term and four-term   expansions 
	  are   an improvement over  the two-term expansion    for $u_c\ll 10^{-4}$  
	 and   
	  a reasonable approximation to $v^*(u_c)$  with acceptable accuracy, less than $10$ percentage error, 
        for $u_c\lesssim 8\times 10^{-3}$ and for $u_c\lesssim 4\times 10^{-3}$, respectively.
  	 However, the higher-order logarithmic corrections that they  neglect are significant for larger values of $u_c$.

\section{Conclusions}

In this paper we have considered a canonical evolution problem for a reaction-diffusion process when the reaction function is of standard KPP-type, but experiences a cut-off in the reaction rate below the normalised cut-off concentration $u_c \in (0,1)$. We have formulated this evolution problem in terms of the moving boundary initial-boundary value problem QIVP. In Section 2 we have obtained some very general results concerning the solution to QIVP. In particular, these general results indicate that in the large-time, as $t \to \infty$, the solution to QIVP will involve the propagation of an advancing non-negative
permanent form travelling wave,
effecting the transition from the unreacting state $u=0$ (ahead of the wave-front) to the fully reacted state $u=1$ (at the rear of the wave-front). With this in mind, this paper has concentrated on examining the existence of permanent form travelling wave  solutions to QIVP with propagation speed $v \geq 0$, referred to as PTW solutions.
In Section 3 we have used a phase plane analysis of the nonlinear boundary value problem \eeref{BVP_cutoff} to establish that (i) for each $u_c \in (0,1)$, then QIVP has a unique PTW solution, with propagation speed $v=v^*(u_c)>0$ and (ii)
$v^*: (0,1) \to \mathbb{R}^+$ is continuous and monotone decreasing, with $v^*(u_c) \to 0^+$ as $u_c \to 1^-$, and $v^*(u_c) \to 2^-$ as $u_c \to 0^+$. It should be noted that $2$ is the minimum propagation speed of permanent form travelling wave solutions for the related KPP-type function in the absence of cut-off.
In Section 4, we have developed asymptotic methods to determine the asymptotic forms of $v^*(u_c)$ as $u_c \to 0^+$ and $u_c \to 1^-$.
 The first limit was previously considered by Brunet and Derrida \cite{BrunetDerrida1997} and Dumortier, Popovic and Kaper \cite{Dumortier_etal2007}. The latter employed
 {geometric desingularisation}
 to  {systematically} determine the order of the error
in \cite{BrunetDerrida1997}.
We have here used matched asymptotics expansions on the direct  problem \eeref{BVP_cutoff} to obtain higher order corrections in a systematic manner. We show that these are controlled by the detailed
structure ahead of the wave-front solution  travelling
with   speed $2$
for the related KPP problem obtained in the absence of a cut-off.
The second limit of $u_c \to 1^-$ is motivated by applications in 
  combustion \cite{Williams1985}.
In this limit, the asymptotic behaviour is obtained via the use of
 regular asymptotic expansions in the  phase plane.

  We anticipate that the approach developed in this paper, for considering PTW solutions to QIVP, will be readily adaptable to corresponding problems, when  the cut-off KPP-type reaction considered here is replaced by a broader class of cut-off   reaction functions, 
  {such as those considered in   \cite{Dumortier_etal2007, Gordon2007,Dumortier_etal2010,Popovic2011,DumortierKaper2012}}. 
  In comparing the PTW theory for the cut-off KPP-type reaction function studied here, and its associated KPP-type reaction function without cut-off, we make the observation that, in the absence of cut-off, a PTW solution exists for each propagation speed $v \in [2, \infty)$, whilst at each fixed cut-off value $u_c \in (0,1)$, a PTW solution exists only at the single propagation speed $v=v^*(u_c)$, with $0 < v^*(u_c)< 2$;  {this observation has been made previously in \cite{Dumortier_etal2007}, although
  restricted to sufficiently small cut-off values $u_c$}. 
   This will have implications for the development of PTW solutions as large-$t$ structures in QIVP, with more general classes of initial data.
In the companion paper we consider the evolution problem QIVP in more detail. Specifically we establish that, as $t \to \infty$, the solution to QIVP does indeed involve the formation of the PTW solution considered in this paper, and we give the detailed asymptotic structure of the solution to QIVP as $t \to \infty$.

Finally, it is interesting to contrast our results with
 results obtained for a related problem,     the stochastic KPP equation 
\begin{subequations}
 \begin{linenomath}\eelab{stochasticKPP}
	\begin{align}
	&	u_t  = u_{xx} + f(u)+(\hat u_c f(u))^{1/2}  \dot W, \quad (x,t) \in \mathbb{R} \times \mathbb{R}^+,\eelab{KPPrand}\\
   & u(x,0)=u_0(x),
 \end{align}
 \end{linenomath}
 \end{subequations}
where $\dot W$ is a standard space-time white-noise.
Similarly to the cut-off KPP equation \eeref{cut-offKPP}, equation \eeref{stochasticKPP} arises as a continuum 
 approximation to (microscopic) interacting particle systems.
 In particular,
for  a Fisher reaction function \eeref{Fisher}, there is an exact relationship between this problem and
 discrete  systems of particles which undergo a birth-coagulation type of reaction
in addition to 
diffusion  
 \cite{ShigaUchiyama1986,Doering_etal2003}. 
Rigorous results have  been derived for this model too  \cite{ConlonDoering2005,Mueller_etal2011}, establishing that the average speed of the random travelling wave solutions of  \eeref{stochasticKPP} is, in the small-$\hat u_c$ or, weak noise limit,  given by
 \beq\eelab{stochasticKPPspeedsmall}
   v_s(\hat u_c)=2- \frac{\pi^2}{(\ln \hat u_c)^2}+O\left(\frac{\ln |\ln \hat u_c|}{|\ln \hat u_c|^3}
  \right),
  \quad \mbox{as} \quad \hat u_c\to 0^+.
  \eeq
  Thus, taking $\hat u_c=u_c$,
  the  difference between  \eeref{stochasticKPPspeedsmall} and the speed
  of the PTW solution of the cut-off KPP model  \eeref{cut-offKPP}
  only arises in the third term of the asymptotic expansion of $v_s(\hat u_c)$ and  $v^*(u_c)$
  as  $u_c \to 0^+$,
a conjecture that was initially made  by
    Brunet and Derrida \cite{BrunetDerrida1997,BrunetDerrida2001}.
 The two models behave very differently 
  when $\hat u_c$    can no longer be regarded as small,  as might be anticipated.
In the large-$\hat u_c$ or, strong noise limit, \cite{Doering_etal2003,HallatschekKorolev2009}
  find that
  \beq\eelab{stochasticKPPspeed}
     v_s(\hat u_c)\sim \sqrt\frac{1}{\hat u_c},
    \quad \mbox{as} \quad \hat u_c\to \infty.
\eeq
 The behaviour in this limit should be contrasted
 against expression \eeref{speed_large_cutoff} obtained for $u_c\to 1^{-}$.
  A comparison suggests that $\hat u_c$ and $u_c$
  may in this case be related according to $\hat u_c\sim 1/(1-u_c)^2$ as $u_c\to 1^-$.
  It would be interesting to extend this comparison to arbitrary $u_c$.

\section*{Acknowledgments}
The research of A D O Tisbury was supported by an EPRSC grant with reference number 1537790. A Tzella thanks C Doering and J Vanneste for useful conversations. 
 {All authors thank the referees and J M Meyer for constructive comments}.

\bibliographystyle{abbrv}
\providecommand{\noopsort}[1]{}\providecommand{\singleletter}[1]{#1}%

\end{document}